\documentclass[12pt]{amsart}
\newcommand{\PP}{{\bf P}}
\newcommand{\RP}{{\bf RP}}
\newcommand{\CP}{{\bf CP}}
\newcommand{\R}{{\bf R}}
\newcommand{\C}{{\bf C}}
\newcommand{\Z}{{\bf Z}}

\newcommand{\sm}{{\setminus}}
\author{V.~A.~Vassiliev}

\thanks{supported in part by RFBR (project 98-01-00555) and INTAS
(grant 96-0713)}

\title{How to calculate homology groups of spaces of
nonsingular algebraic projective hypersurfaces}

\date{Revised version published in 1999}

\begin{document}

\begin{abstract}
A general method of computing cohomology groups of the space of nonsingular
algebraic hypersurfaces of degree $d$ in $\CP^n$ is described. Using this
method, rational cohomology groups of such spaces with $n=2, d \le 4$  and
$n=3=d$ are calculated.
\end{abstract}

\maketitle

\section{Introduction}

I describe a method of calculating the cohomology groups of spaces of
nonsingular algebraic hypersurfaces of fixed degree in $\CP^n$: a complex
analog of the rigid homotopy classification of algebraic hypersurfaces in
$\RP^n.$

This method is essentially the one described in \cite{book}, \cite{phasis}; it
consists of two main parts. The first one is due to Arnold, who remarked
\cite{ar70} that often it is convenient to replace the calculation of
cohomology groups of spaces of nonsingular objects by that of the (Alexander
dual) homology groups of corresponding {\em discriminant} spaces of singular
objects; see also \cite{ar89}. The second one is the techniques of conical
resolutions and topological order complexes, which is a continuous analog of
the combinatorial inclusion-exclusion formula and is adopted very well to the
calculation of the latter homology groups. This construction was proposed in
\cite{viniti} and generalized in \cite{v91}, \cite{book}, \cite{phasis}.

The version of this method, applied to spaces of nonsingular projective
hypersurfaces, is described in \S~3. In the preceding \S~2 we demonstrate all
main ingredients of it, considering the simplest problem of this series: we
calculate there the (well-known) cohomology groups of the space of regular
quadrics in $\CP^2.$ In \S~4 we do the same for the real cohomology groups of
the space of nonsingular cubics in $\CP^2$ (which probably also can be
calculated by more standard methods, see remark after Theorem 2). In \S~5, \S~6
and \S~7 we present similar calculations respectively for spaces of nonsingular
quartics in $\CP^2$, nonsingular cubics in $\CP^3$, and nondegenerate quadratic
vector fields in $\C^3$.
\medskip

This work is inspired by the Arnold's problem 1970-13 on the topology of the
space of nonsingular cubics, see \cite{aprobl}.

\subsection{Stating the problem and the first reduction.}

Let $\Pi_{d,n} \simeq \C^{D(d,n)}$ be the space of all homogeneous polynomials
$\C^{n+1} \to \C^1$ of degree $d,$ $\Sigma_{d,n} \subset \Pi_{d,n}$ the set of
forms, having singular points outside the origin, and $N_{d,n} \subset
\CP^{D(d,n)-1}$  the projectivization of $\Pi_{d,n} \setminus \Sigma_{d,n}$,
i.e. the set of nonsingular algebraic hypersurfaces of degree $d$ in $\CP^{n}.$
$\Sigma_{d,n}$  is a conical nonempty algebraic hypersurface in $\Pi_{d,n},$
therefore $\Pi_{d,n} \setminus \Sigma_{d,n}$ is diffeomorphic to $N_{d,n}
\times \C^*$, in particular the homology groups of spaces $N_{d,n}$ and
$\Pi_{d,n} \setminus \Sigma_{d,n}$ are related via the K\"unneth formula.
\medskip

By the Alexander duality theorem, for $i>0$  we have
\begin{equation}
\label{alex} \tilde H^i(\Pi_{d,n} \setminus \Sigma_{d,n}) \simeq \bar
H_{2D(d,n) -i-1}(\Sigma_{d,n}),
\end{equation}
where $\bar H_*$ is the notation for the Borel--Moore homology, i.e. the
homology of the one-point compactification reduced modulo the added point.

\section{The trivial example: nonsingular quadrics in $\CP^2$}

In this section we assume that $n=2$ and denote the spaces $\Pi_{d,2},
\Sigma_{d,2}$ and $N_{d,2}$ simply by $\Pi_{d}, \Sigma_{d}$ and $N_{d}$.
\medskip

{\sc Proposition 1.} {\em Only the following groups $H^i(\Pi_2\setminus
\Sigma_2)$ are nontrivial:} $H^0 \simeq H^1 \simeq H^5 \simeq H^6 \simeq \Z,$
$H^3 \simeq H^4 \simeq \Z_2.$ {\em In particular, only the following groups
$H^i(N_2)$ are nontrivial:} $H^0 \simeq H^5 \simeq \Z,$ $H^3 \simeq \Z_2.$
\medskip

This statement is not new: indeed, it is easy to see that the space $\Pi_2
\setminus \Sigma_2$ is homotopy equivalent to the Lagrangian Grassmannian
manifold $U(3)/O(3),$ whose cohomology groups were studied in \cite{borel} and
\cite{fuchs2}, see e.g. \cite{v88}.  We give here another calculation,
demonstrating our general method.

This calculation is based on the classification of subsets in $\CP^2$, which
can be the singular sets of homogeneous polynomials of degree $2$ in $\C^3.$
There are exactly the following such sets:

A) any point $x \in \CP^2$,

B) any line $\CP^1 \subset \CP^2$ (such lines are parametrized by the points of
the dual projective space $\CP^{2\vee}$);

C) the entire $\CP^2$.

For any set of type A), B) or C), the corresponding set of quadratic forms
$\C^3 \to \C$, having singularities at all points of this set (and maybe
somewhere else) is isomorphic, respectively, to $\C^3$, $\C^1$ or $\C^0.$

The main topological tool of the calculation is the following one.
\medskip

\subsection{The topological order complex $\Lambda_2$}

Take the spaces $\CP^2, \CP^{2\vee}$ and $\{\bullet\}$  of sets of types A) B)
and C), and consider the join $\CP^2 * \CP^{2\vee} * \{\bullet\} \ $ of these
spaces, i.e., roughly speaking, the union of all simplices (of dimensions 0, 1
and 2), whose vertices correspond to their points.

Such a simplex is called {\em coherent,} if all subsets in $\CP^2$,
corresponding to its vertices, are incident to one another.

The {\em topological order complex} $\Lambda_2$ is defined as the union of all
coherent simplices, with topology induced from that of the join.

To any set $K$ of the form A), B) or C) there corresponds a subset
$\Lambda_2(K) \subset \Lambda_2:$ the union of simplices, all whose vertices
are subordinate by inclusion to the point $\{K\}$.  Obviously this subset is a
cone with the vertex $\{K\}$:  for $K$ of type A) (respectively, B),
respectively, C)) it is equal to $\{K\}$ itself (respectively, to the cone over
the topological space $K \simeq \CP^1$, respectively, to entire $\Lambda_2$).
By ${\stackrel{\circ}{\Lambda}}_2(K)$ we denote this cone with its base
removed.

\subsection{Conical resolution}
The {\em conical resolution} $\sigma_2$ of $\Sigma_2$ is defined as a subset in
the direct product $\Pi_2 \times \Lambda_2$.  Namely, let $K\subset \CP^2$ be
any possible set of type A), B) or C). Let $L(K)$ be the set of all quadratic
forms $\C^3 \to \C$, having singular points at all points of $K$ (and maybe
somewhere else); this always is a vector subspace in $\Pi_2.$ Then the space
$\sigma_2(K) \subset \Pi_2 \times \Lambda_2$ is defined as the direct product
$L(K) \times \Lambda_2(K),$ and the desired conical resolution $\sigma_2$ as
the union of all possible sets $\sigma_2(K).$
\medskip

{\sc Proposition 2} (see \cite{book}, \cite{phasis}). {\em The obvious
projection $\sigma_2 \to \Sigma_2$ is a proper map, and the induced map of
one-point compactifications of these spaces  is a homotopy equivalence. In
particular this projection induces the isomorphism $\bar H_*(\sigma_2) \simeq
\bar H_*(\Sigma_2)$.} \quad $\Box$ \medskip

The space $\sigma_2$ admits a natural filtration $F_1 \subset F_2 \subset F_3$:
the subspace $F_1$ (respectively, $F_2$, respectively, $F_3$) of $\sigma_2$ is
the union of sets $\sigma_2(K)$ such that the corresponding spaces $K$ are of
type A) only (respectively, of type A) or B), respectively, of type A) or B) or
C), i.e. $F_3 = \sigma_2$).

Similarly, we consider the filtration $\{\Phi_i\}$ of the space $\Lambda_2:$
the subspace $\Phi_1$ (respectively, $\Phi_2$, respectively, $\Phi_3$) of
$\Lambda_2$ is the union of sets $\Lambda_2(K)$ such that the corresponding
spaces $K$ are of type A) only (respectively, of type A) or B), respectively,
of type A) or B) or C), i.e. $\Phi_3 = \sigma_2$).

In particular, the space $ F_i \setminus F_{i-1}$ is the space of a complex
vector bundle over $ \Phi_i \sm \Phi_{i-1},$ the dimension of this bundle is
equal to $3$ for $i=1,$ to $1$ for $i=2$ and to $0$  for $i=3$. Thus the
Borel--Moore homology groups of spaces $ F_i \setminus F_{i-1}$ and $ \Phi_i
\sm \Phi_{i-1}$ are related by the Thom isomorphism.

Consider the spectral sequence generated by the filtration $\{F_i\}$ and
converging to the group $\bar H_*(\sigma_2).$ By definition its term
$E^1_{p,q}$ is equal to $\bar H_{p+q}(F_p \setminus F_{p-1}).$
\medskip

\begin{figure}
\unitlength=1.00mm \thinlines \linethickness{0.4pt}
\begin{picture}(105.00,60.00)
\put(10.00,10.00){\vector(0,1){50.00}} \put(10.00,10.00){\vector(1,0){40.00}}
\put(6.00,59.00){\makebox(0,0)[cc]{$q$}}
\put(6.00,52.00){\makebox(0,0)[cc]{$9$}}
\put(6.00,47.00){\makebox(0,0)[cc]{$8$}}
\put(6.00,42.00){\makebox(0,0)[cc]{$7$}}
\put(6.00,37.00){\makebox(0,0)[cc]{$6$}}
\put(6.00,32.00){\makebox(0,0)[cc]{$5$}}
\put(6.00,27.00){\makebox(0,0)[cc]{$4$}}
\put(6.00,22.00){\makebox(0,0)[cc]{$3$}}
\put(15.00,32.00){\makebox(0,0)[cc]{${\bf Z}$}}
\put(15.00,42.00){\makebox(0,0)[cc]{${\bf Z}$}}
\put(15.00,52.00){\makebox(0,0)[cc]{${\bf Z}$}}
\put(25.00,42.00){\makebox(0,0)[cc]{${\bf Z}$}}
\put(25.00,32.00){\makebox(0,0)[cc]{${\bf Z}$}}
\put(25.00,22.00){\makebox(0,0)[cc]{${\bf Z}$}}
\put(35.00,32.00){\makebox(0,0)[cc]{${\bf Z}$}}
\put(30.00,3.00){\line(0,1){54.00}} \put(20.00,57.00){\line(0,-1){54.00}}
\put(15.00,6.00){\makebox(0,0)[cc]{$1$}}
\put(25.00,6.00){\makebox(0,0)[cc]{$2$}}
\put(35.00,6.00){\makebox(0,0)[cc]{$3$}}
\put(48.00,6.00){\makebox(0,0)[cc]{$p$}} \put(65.00,10.00){\vector(0,1){50.00}}
\put(65.00,10.00){\vector(1,0){40.00}}
\put(61.00,59.00){\makebox(0,0)[cc]{$q$}}
\put(61.00,52.00){\makebox(0,0)[cc]{$9$}}
\put(61.00,47.00){\makebox(0,0)[cc]{$8$}}
\put(61.00,42.00){\makebox(0,0)[cc]{$7$}}
\put(61.00,37.00){\makebox(0,0)[cc]{$6$}}
\put(61.00,32.00){\makebox(0,0)[cc]{$5$}}
\put(61.00,27.00){\makebox(0,0)[cc]{$4$}}
\put(61.00,22.00){\makebox(0,0)[cc]{$3$}}
\put(70.00,32.00){\makebox(0,0)[cc]{${\bf Z}$}}
\put(70.00,52.00){\makebox(0,0)[cc]{${\bf Z}$}}
\put(80.00,22.00){\makebox(0,0)[cc]{${\bf Z}$}}
\put(85.00,3.00){\line(0,1){54.00}} \put(75.00,57.00){\line(0,-1){54.00}}
\put(70.00,6.00){\makebox(0,0)[cc]{$1$}}
\put(80.00,6.00){\makebox(0,0)[cc]{$2$}}
\put(90.00,6.00){\makebox(0,0)[cc]{$3$}}
\put(103.00,6.00){\makebox(0,0)[cc]{$p$}}
\put(70.00,42.00){\makebox(0,0)[cc]{${\bf Z}_2$}}
\put(80.00,32.00){\makebox(0,0)[cc]{${\bf Z}_2$}}
\end{picture}
\caption{Terms $E^1$ and $E^2$ of the spectral sequence for $d=n=2$.}
\label{ss2}
\end{figure}
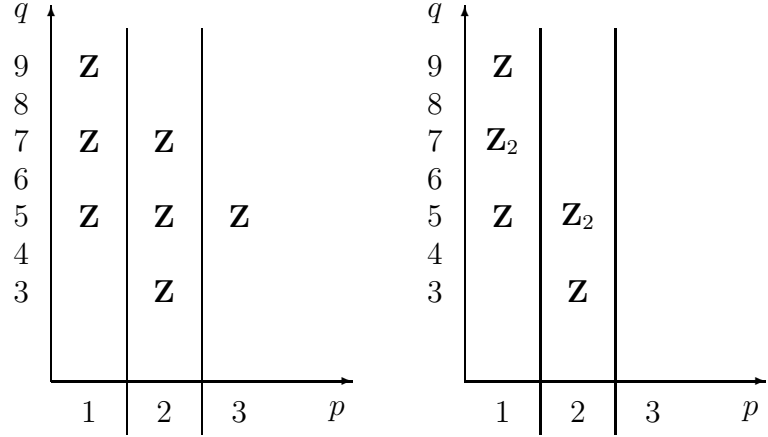

{\sc Lemma 1.} {\em All nontrivial terms $E^1_{p,q}$ of this spectral sequence
are as follows (see fig.~\ref{ss2} left): $E^1_{1,q}=\Z$ for $q=5,7,9;$
$E^1_{2,q}=\Z$ for $q=3,5,7;$ $E^1_{3,q}=\Z$ for $q=5.$}
\medskip

{\em Proof.} 1. $F_1$ is the space of a fiber bundle over $\CP^2$ with fiber
$\C^3,$ this proves the statement concerning $E^1_{1,q}.$

2. $F_2 \setminus F_1$ is fibered over the space $\CP^{2\vee}$ of all complex
lines in $\CP^2.$ The fiber over such a line $K$ is the direct product of the
space $L(K) \simeq \C^1$ and the space $\Lambda_2(K) \setminus \Phi_1.$ But
$\Lambda_2(K)$ is the cone over the set of all points of $K \sim \CP^1,$ and
the base of this cone belongs to $\Phi_1.$ Therefore $\Lambda_2(K) \setminus
\Phi_1$ is an oriented 3-dimensional open disc, thus by the Thom isomorphism
$\bar H_*(F_2 \sm F_1) \simeq \bar H_{*-5}(\CP^{2\vee}),$ and the statement
concerning $E^1_{2,*}$ is proved.

3. $F_3 \setminus F_2 \equiv \Phi_3 \setminus \Phi_2.$ $\Phi_2$ is the subset
of the join $\CP^2 * \CP^{2\vee},$ consisting of all segments, spanning such
pairs $(x \in \CP^2,K \in \CP^{2\vee}),$ that $x \in K.$ In \cite{v91} it was
proved that this subset is homeomorphic to $S^7.$ $\Phi_3$ is the cone over
$\Phi_2,$ hence $F_3 \setminus F_2$ is an open 8-dimensional disc, and the
statement concerning $E^1_{3,*}$ also is proved.
\medskip

{\sc Proposition 3.} {\em The differential $d^1$ of our spectral sequence maps
a generator of the group $E^1_{3,5}$ into twice a generator of $E^1_{2,5},$
maps a generator of $E^1_{2,7}$ into twice a generator of $E^1_{1,7},$ and acts
trivially on all other groups $E^1_{p,q}$. In particular all nontrivial groups
$E^2_{p,q}$ are as shown in fig.~\ref{ss2} right and our spectral sequence
degenerates at this step: $E^2 \equiv E^\infty$.}
\medskip

Proposition 1 follows immediately from this one, Lemma 1 and the Alexander
isomorphism (\ref{alex}).

\medskip
{\em Proof of Proposition 3.} 1. The space $F_2 \sm F_1$ can be considered as
the space of a fiber bundle with fiber $\C^1,$ whose base $\Phi_2 \sm \Phi_1$
is fibered over $\CP^{2\vee}$ with fiber equal to an open 3-dimensional disc.
The image of the map $d^1:E^1_{3,5} \to E^1_{2,5}$ is equal to the fundamental
cycle (with closed supports) of the base $\Phi_2 \setminus \Phi_1$  of this
bundle. The homomorphism
\begin{equation}
\label{thom1} \bar H_7(\Phi_2 \sm \Phi_1) \to \bar H_7(F_2 \sm F_1) \to \bar
H_5(\Phi_2 \sm \Phi_1),
\end{equation}
defined as the composition of the identical embedding and the Thom isomorphism,
maps this fundamental cycle into a homology class, Poincar\'e dual to the Euler
class of the complex line bundle, whose fiber over the points of the open cone
${\stackrel{\circ}{\Lambda}}_2(K),$ $K \in \CP^{2\vee},$ is equal to the space
$L(K).$ This space consists of all quadratic functions of the form $l^2,$ where
$l$ is an arbitrary linear form vanishing on the line $K.$ It is easy to
calculate that this Euler class is equal to twice the generator of the group
$\bar H_5(\Phi_2 \sm \Phi_1),$ and the first statement of Proposition 3 is
proved.

2. The homomorphism $d^1:E^1_{2,7} \to E^1_{1,7}$ maps the generator of the
group $E^1_{2,7}$ into the fundamental cycle of the algebraic subset $\omega
\subset F_1$, consisting of all pairs of the form $($a point $u \in \CP^2$, a
quadratic function $l^2),$ such that $l$ is a linear form with $l(u)=0.$ This
subset is fibered over $\CP^2$ with a standard quadratic cone $\{\alpha \beta =
\gamma^2\} \subset \C^3$  for a fiber. To calculate the homology class of this
cycle, it is sufficient to calculate its intersection index with the generator
of the group $H_2(F_1).$ This generator can be realized as follows: we take an
arbitrary line $\CP^1 \subset \CP^2,$ say the one given by the equation $x=0,$
and consider any section of the restriction on it of the fibration $F_1 \to
\CP^2.$ For such a section we can take the one, whose value over the point
$(0:1:\lambda) \in \CP^1$ is the quadratic form $\frac{x(z-\lambda
y)}{1+|\lambda|^2}.$ This section intersects the above-described set $\omega$
at exactly one point, corresponding to $\lambda=\infty.$ This intersection is
of multiplicity 2, hence the fundamental cycle of this set $\omega$ is equal to
twice the generator of $\bar H_8(F_1)$.

It is now obvious that $d^1$ cannot be nontrivial on any other group
$E^1_{p,q},$ $(p,q) \ne (3,5), (2,7),$ and Proposition 3 is completely proved.

\section{The general spectral sequence}

Here we extend the above-described construction to the case of hypersurfaces of
an arbitrary degree $d$. The spaces of sets of types A), B) and C) will be
formalized in terms of the Hilbert-scheme topology.

\subsection{Conical resolution of the discriminant:
preliminary version.}

Set $D\equiv D(d,n) =\dim \Pi_{d,n}$. Consider all Grassmann manifolds
$G_i(\Pi_{d,n}),$ $i=1, \ldots, D,$ whose points are the complex subspaces of
{\em codimension} $i$ in $\Pi_{d,n}$. Suppose that the subset $K \subset \CP^n$
is the set of singular points of a certain homogeneous function $\C^{n+1} \to
\C^1,$ and denote by $L(K)$  the subspace in $\Pi_{d,n}$ consisting of all
functions having singular points everywhere in $K$ (and maybe somewhere else).
Obviously this is a vector subspace in $\Pi_{d,n}$. Let $\Omega_i \subset
G_i(\Pi_{d,n})$ be the set of all subspaces of the form $L(K).$ Let $\bar
\Omega_i$ be the closure of $\Omega_i$ in $G_i(\Pi_{d,n})$.

Consider the join $G_1(\Pi_{d,n}) * \cdots
* G_{D-1}(\Pi_{d,n}) * G_D(\Pi_{d,n})$
and the subset in it, consisting of all {\em coherent} simplices of any
dimensions, i.e. of such simplices, that a) all their vertices belong to
subspaces $\bar \Omega_i$, and b) the collection of subspaces in $\Pi_{d,n}$,
corresponding to all vertices of the simplex, forms a flag.

The union of all such simplices will be denoted by $\tilde \Lambda_{d,n}$. The
{\em leading vertex} of a coherent simplex is the one lying in the space $\bar
\Omega_i$ with the greatest value of $i$. For any point $\{L\} \in \bar
\Omega_i,$ the set $\tilde \Lambda_{d,n}(L)$ is the union of coherent
simplices, whose leading vertices coincide with $L$. Obviously this is a cone
with vertex $\{L\}$. In particular, $\tilde \Lambda_{d,n} = \tilde
\Lambda_{d,n}(\{\Pi_{d,n}\})$.

The subset $\tilde \sigma_{d,n}(L) \subset \Pi_{d,n} \times \tilde
\Lambda_{d,n}$ is defined as the direct product $L \times \tilde
\Lambda_{d,n}(L).$ The {\em conical resolution} $\tilde \sigma_{d,n}$ of
$\Sigma_{d,n}$ is defined as the union of all sets $\tilde \sigma_{d,n}(L)$
over all points $\{L\}$ of all $\bar \Omega_i,$ $i=1, \ldots, D.$
\medskip

{\sc Proposition 4.} {\em The obvious projection $\tilde \sigma_{d,n} \to
\Sigma_{d,n}$ is a proper map, and the induced map of one-point
compactifications of these spaces  is a homotopy equivalence. In particular it
induces an isomorphism $\bar H_*(\tilde \sigma_{d,n}) \simeq \bar
H_*(\Sigma_{d,n})$.}
\medskip

The space $\tilde \sigma_{d,n}$ admits a natural increasing filtration: its
term $\tilde F_i \subset \tilde \sigma_{d,n}$ is the union of sets $\tilde
\sigma_{d,n}(L)$ over all $L \in \bar \Omega_j,$ $j \le i$.

Simultaneously we consider the increasing filtration $\{\tilde \Phi_i\}$ on the
space $\tilde \Lambda_{d,n}:$ its term $\tilde \Phi_i$ is the union of sets
$\tilde \Lambda_{d,n}(L)$ such that $\mbox{codim}_{\C}L \le i.$ In particular,
the space $\tilde F_i \setminus \tilde F_{i-1}$ is the space of a complex
$(D-i)$-dimensional vector bundle over $\tilde \Phi_i \setminus \tilde
\Phi_{i-1},$ and their Borel--Moore homology groups are related by the Thom
isomorphism.

The simplicial resolution of \S~2  is the special case of this one, and its
filtration $F_1 \subset F_2 \subset F_3$ is obtained from our present
filtration $\{\tilde F_i\}$ by renumbering its terms: $F_1=\tilde F_3 =\tilde
F_4,$ $F_2=\tilde F_5,$ $F_3=\tilde F_6.$

\subsection{Reduced conical resolutions of discriminants.}

There is another, more economical conical resolution of $\Sigma_{d,n}$, also
coinciding with the previous one in the case $d=2=n$. Namely, we can reduce all
``nongeometrical'' vertices $L \in \bar \Omega_i \sm \Omega_i$ and their
subordinate simplices.

\medskip
{\sc Definition.} Given a point $L \in \bar \Omega_i,$ its {\em geometrization}
$\bar L$ is the smallest plane of the form $L(K) \in \Omega_j,$ $j \le i,$
containing $L$. In particular $\bar L = L$ if $L \in \Omega_i.$
\medskip

{\sc Example.}  Suppose that $d \ge 3,$ $u \in \CP^n,$ and $v(t),$ $t \in
\C^1,$ is an one-parametric analytic family of points in $\CP^n$ such that
$v(t)=u \Leftrightarrow t=0.$ Then for any $t \ne 0$ the 2-point set $K(t)
\equiv (u,v(t)) $ defines a subspace $L(t) \equiv L(K(t))$ of codimension
$2(n+1)$ in $\Pi_{d,n},$ while the limit position $L(0)$ of these subspaces
when $t \to 0$ belongs to $\bar \Omega_{2(n+1)} \sm \Omega_{2(n+1)}.$ Then
$\overline{L(0)} = L(\{u\}) \in \Omega_{n+1}.$

The segment in $\tilde \Lambda_{d,n}$,  joining these two points $L(0)$ and
$L(\{u\})$, is homotopy equivalent to its endpoint $L(\{u\}),$ therefore we can
contract it to this point, and to extend this retraction by linearity to all
coherent simplices, containing this segment. Doing the same in all similar
situations, we reduce the complex $\tilde \Lambda_{d,n}$ to the (properly
topologized) union $\Lambda_{d,n}$ of all coherent simplices, whose vertices
correspond to the points of spaces $\Omega_i$ only, but not of $\bar \Omega_i
\sm \Omega_i.$ Then we lift this reducing map to the space $\tilde
\sigma_{d,n}$ and get the desired space $\sigma_{d,n},$ with which we will work
all the time. This general reduction map is constructed as follows.

\medskip
Given a point $L \in \Omega_j,$ denote by $W(L)$ the union of all coherent
simplices in $\tilde \Lambda_{d,n}$ such that for all their vertices $L'$ their
geometrizations $\overline{L'}$  coincide with $L$. Obviously $W(L)$ is a
compact and contractible subspace of $\tilde \Lambda_{d,n}$: it is a cone with
vertex $\{L\}.$
\medskip

The {\it reducing map}
\begin{equation}
\label{redd} red: \tilde \Lambda_{d,n} \to \Lambda_{d,n}
\end{equation}
is defined as follows. For any $i$ and any $L \in \bar \Omega_i$ we map the
point $L$ to $\bar L$ and extend this map by linearity to all coherent
simplices; the image of a coherent simplex under this map obviously is again
coherent. Then we call {\em equivalent} all points of $\tilde \Lambda_{d,n}$
sent by this extended map into one and the same point. The space
$\Lambda_{d,n}$ is defined as the quotient space of $\tilde \Lambda_{d,n}$ by
this equivalence relation.

This factorization map defines obviously a map $\tilde \Lambda_{d,n} \times
\Pi_{d,n} \to \Lambda_{d,n} \times \Pi_{d,n},$ and hence  also a map
\begin{equation}
\label{red} Red: \tilde \sigma_{d,n} \to  \sigma_{d,n}
\end{equation}
of the resolved discriminant $\tilde \sigma_{d,n} \subset \tilde \Lambda_{d,n}
\times \Pi_{d,n}$ onto some space $\sigma_{d,n} \subset \Lambda_{d,n} \times
\Pi_{d,n}$.
\medskip

{\sc Proposition 5.} {\it The map $($\ref{red}$)$ induces a homotopy
equivalence of one-point compactifications of spaces $\tilde \sigma_{d,n}$ and
$ \sigma_{d,n}$ and factorizes the map $\pi: \tilde \sigma_{d,n} \to \Sigma$
(i.e., there is a proper map $\tilde \pi:  \sigma_{d,n} \to \Sigma$ such that
$\pi = \tilde \pi \circ Red$).}
\medskip

This follows immediately from the construction. \quad $\Box$
\medskip

For any $L \in \Omega_i$  denote by $\Lambda_{d,n}(L) $ the image of $\tilde
\Lambda_{d,n}(L)$  under the reducing map (\ref{redd}).
\medskip

{\sc Example.} If $K$ is a set of $r < \infty$  points, and $d$ is sufficiently
large with respect to $r,$ then  both spaces $\tilde \Lambda_{d,n}(L(K))$ and $
\Lambda_{d,n}(L(K))$ are naturally homeomorphic to a $(r-1)$-dimensional
simplex, and the restriction of the map $red$ onto $\tilde \Lambda_{d,n}(K(L))$
is a homeomorphism between them.
\medskip

{\sc Definitions and notation.} For any set $L$ we denote by $\partial
\Lambda_{d,n}(L)$ the {\em link}  of the order complex $\Lambda_{d,n}(L),$ i.e.
the union of all its coherent simplices not containing its vertex $\{L\}$, and
denote by $\stackrel{\circ}{\Lambda}_{d,n}(L)$ the {\em open cone}
$\Lambda_{d,n}(L) \sm \partial \Lambda_{d,n}(L).$ Their homology groups are
obviously related by the boundary isomorphism
\begin{equation}
\label{bound} \bar H_*(\stackrel{\circ}{\Lambda}_{d,n}(L)) \cong
H_*(\Lambda_{d,n}(L),
\partial \Lambda_{d,n}(L)) \stackrel{\sim}{\to}
H_{*-1}(\partial \Lambda_{d,n}(L),pt).
\end{equation}
For any subset $K \subset \CP^n$ the complex $\Lambda_{d,n}(L(K))$ will be
usually denoted simply by $\Lambda_{d,n}(K)$.

\subsection{Filtrations and spectral sequences.}

The space $\tilde \Lambda_{d,n}$ admits a natural increasing filtration $\tilde
\Phi_1 \subset \cdots \subset \tilde \Phi_D \equiv \tilde \Lambda_{d,n}$.
Namely, the filtration of any interior point of a coherent simplex $\Delta
\subset \tilde \Lambda_{d,n}$ is equal to the codimension $i$ of the plane
$L(\Delta) \in \bar \Omega_i$ corresponding to the leading vertex of this
simplex. Hence, a similar filtration is induced also on the quotient space
$\Lambda_{d,n}$: the filtration of the equivalence class $x \in \Lambda_{d,n}$
is equal to the smallest value of filtrations of points in $\tilde
\Lambda_{d,n}$, constituting $x$.  The obvious projections $\tilde \sigma_{d,n}
\to \tilde \Lambda_{d,n}$, $ \sigma_{d,n} \to \Lambda_{d,n}$ induce also
filtrations on spaces $\tilde \sigma_{d,n}$ and $ \sigma_{d,n}$. By
construction, the term $\tilde F_i \setminus \tilde F_{i-1}$ of this filtration
on $\tilde \sigma_{d,n}$ (respectively, on $ \sigma_{d,n}$) is the space of a
$(D-i)$-dimensional complex vector bundle over the term $\tilde \Phi_i \sm
\tilde \Phi_{i-1}$ of the filtration on $\tilde \Lambda_{d,n}$ (respectively,
on $\Lambda_{d,n}$).

In fact, we will use a shortened version of this filtration, removing all
non-increasing terms by the following inductive rule: its term $F_1$ is the
first nonempty set $\tilde F_j$, and its term $F_i$ is the first term $\tilde
F_j$ of the former filtration which is {\em strictly} greater than $F_{i-1}$.

This filtration defines a spectral sequence, converging to the right-hand group
in (\ref{alex}). Its term $E^1_{p,q}$ is isomorphic to $\bar H_{p+q}(F_p
\setminus F_{p-1})$.

\section{Non-singular forms of degree 3}

{\sc Theorem 1.} {\em Only the following groups $H^i(\Pi_{3,2}\setminus
\Sigma_{3,2},\R)$ are nontrivial:} $H^0 \simeq H^1 \simeq H^3 \simeq H^4 \simeq
H^5 \simeq H^6 \simeq H^8 \simeq H^9 \simeq \R.$ {\em In particular, only the
following groups $H^i(N_{3,2},\R)$ are nontrivial:} $H^0 \simeq H^3 \simeq H^5
\simeq H^8 \simeq \R.$
\medskip

This theorem follows immediately from the following one.

\medskip
{\sc Theorem 2.} {\em The term $E^1$ of our spectral sequence (over $\R$),
converging to the group $\bar H_*(\Sigma_{3,2}),$ looks as in fig.~\ref{ss3}
(where empty cells denote trivial groups $E^1_{p,q}).$ $E^\infty \equiv E^1$,
i.e. the unique possible differential $d_1:E^1_{2,13} \to E^1_{1,13}$ of this
sequence is trivial.}
\medskip

\begin{figure}
\unitlength=1.00mm \thinlines \linethickness{0.4pt}
\begin{picture}(95.00,100.00)
\put(10.00,10.00){\vector(1,0){85.00}} \put(93.00,6.00){\makebox(0,0)[cc]{$p$}}
\put(10.00,10.00){\vector(0,1){90.00}} \put(5.00,98.00){\makebox(0,0)[cc]{$q$}}
\put(5.00,87.00){\makebox(0,0)[cc]{$17$}}
\put(5.00,81.00){\makebox(0,0)[cc]{$16$}}
\put(5.00,75.00){\makebox(0,0)[cc]{$15$}}
\put(5.00,69.00){\makebox(0,0)[cc]{$14$}}
\put(5.00,63.00){\makebox(0,0)[cc]{$13$}}
\put(5.00,57.00){\makebox(0,0)[cc]{$12$}}
\put(5.00,51.00){\makebox(0,0)[cc]{$11$}}
\put(5.00,45.00){\makebox(0,0)[cc]{$10$}}
\put(5.00,39.00){\makebox(0,0)[cc]{$9$}}
\put(5.00,33.00){\makebox(0,0)[cc]{$8$}}
\put(5.00,27.00){\makebox(0,0)[cc]{$7$}}
\put(5.00,21.00){\makebox(0,0)[cc]{$6$}} \put(22.00,4.00){\line(0,1){89.00}}
\put(34.00,93.00){\line(0,-1){89.00}} \put(46.00,4.00){\line(0,1){89.00}}
\put(58.00,93.00){\line(0,-1){89.00}} \put(70.00,4.00){\line(0,1){89.00}}
\put(16.00,87.00){\makebox(0,0)[cc]{${\bf R}$}}
\put(16.00,75.00){\makebox(0,0)[cc]{${\bf R}$}}
\put(16.00,63.00){\makebox(0,0)[cc]{${\bf R}$}}
\put(28.00,63.00){\makebox(0,0)[cc]{${\bf R}$}}
\put(28.00,51.00){\makebox(0,0)[cc]{${\bf R}$}}
\put(28.00,39.00){\makebox(0,0)[cc]{${\bf R}$}}
\put(52.00,21.00){\makebox(0,0)[cc]{${\bf R}$}}
\put(64.00,6.00){\makebox(0,0)[cc]{$5_{(10)}$}}
\put(52.00,6.00){\makebox(0,0)[cc]{$4_{(9)}$}}
\put(40.00,6.00){\makebox(0,0)[cc]{$3_{(7)}$}}
\put(28.00,6.00){\makebox(0,0)[cc]{$2_{(6)}$}}
\put(16.00,6.00){\makebox(0,0)[cc]{$1_{(3)}$}}
\put(2.00,18.00){\line(1,0){69.00}} \put(2.00,24.00){\line(1,0){69.00}}
\put(2.00,30.00){\line(1,0){69.00}} \put(2.00,36.00){\line(1,0){69.00}}
\put(2.00,42.00){\line(1,0){69.00}} \put(2.00,48.00){\line(1,0){69.00}}
\put(2.00,54.00){\line(1,0){69.00}} \put(2.00,60.00){\line(1,0){69.00}}
\put(2.00,66.00){\line(1,0){69.00}} \put(2.00,72.00){\line(1,0){69.00}}
\put(2.00,78.00){\line(1,0){69.00}} \put(2.00,84.00){\line(1,0){69.00}}
\put(2.00,90.00){\line(1,0){69.00}}
\end{picture}
\caption{spectral sequence for $d=3$.} \label{ss3}
\end{figure}

{\sc Remark.} Theorem 1 is not surprising: it states in fact that
$H^*(\Pi_{3,2} \sm \Sigma_{3,2}, \R) \simeq H^*(PGL(\CP^2),\R).$ There is
obvious projection of the space $\Pi_{3,2} \sm \Sigma_{3,2}$ onto the modular
space of elliptic curves. The group $PGL(\CP^2)$ acts transitively on all
fibers of this projection, and the stationary groups of this action are finite
for all fibers (although they are not the same: they "jump" over the curve with
additional symmetry). It seems likely that Theorem 1 could be proved by
considering these subgroups and the Leray spectral sequence of this projection.
\medskip

The proof of Theorem 2 occupies the rest of this section. It uses many times
the following notions and facts.

\subsection{Main homological lemmas.}

{\sc Definition.} Given a topological space $X$, the {\em $k$-th configuration
space $B(X,k)$ of $X$} is the space of all subsets of cardinality $k$ in $X,$
supplied with the natural topology, see e.g. \cite{book}, \cite{phasis}. The
{\em sign representation} $\pi_1(B(X,k)) \to \mbox{Aut} (\Z)$ maps the paths in
$B(X,k),$  defining odd (respectively, even) permutations of $k$ points into
multiplication by $-1$ (respectively, $1$). The local system $\pm \Z$ over
$B(X,k),$ is the one locally isomorphic to $\Z$, but with monodromy
representation equal to the sign representation of $\pi_1(B(X,k))$. The {\em
twisted Borel--Moore homology group} $\bar H_*(B(X,k),\pm \Z)$ is defined as
the homology group of the complex of {\em locally finite} singular chains of
$B(X,k)$ with coefficients in the sheaf $\pm \Z$. Similarly we define local
systems $\pm \R \equiv \pm \Z \otimes \R$ and their homology groups. For any
topological space $X,$ $\bar X$ denotes the one-point compactification of $X$.
\medskip

{\sc Lemma 2.} {\em A. The group
\begin{equation}
\label{cnk} \bar H_*(B(\C^n,k), \pm \R),
\end{equation}
$n \ge 1,$ is trivial for any $k \ge 2.$

B. The group
\begin{equation}
\label{cpnk} \bar H_*(B(\CP^n,k), \pm \R),
\end{equation}
$n \ge 1,$ is isomorphic to $H_{*-k(k-1)}(G_k(\C^{n+1}),\R),$ where
$G_k(\C^{n+1})$ is the Grassmann manifold of $k$-dimensional subspaces in
$\C^{n+1},$ see} \cite{MS}.

{\em In particular, the group (\ref{cpnk}) is trivial if $k>n+1$.}
\medskip

{\em Proof.} Assertion A is well-known, see e.g. \cite{book}, \cite{phasis}.

Let us fix some complete flag $\bullet \subset \CP^1 \subset \CP^2 \subset
\cdots \subset \CP^n.$ To any configuration of $k$ points in $\CP^n$ associate
its {\em index} $a=(a_0 \le a_1 \le \cdots \le a_n),$ $a_n=k,$ where $a_i$ is
the number of its points lying in $\CP^i.$ For any such index $a$ denote by
$B_a \subset B(\CP^n,k)$ the union of all configurations having this index. It
follows immediately from statement A, that $\bar H_*(B_a,\pm \R) =0$ for any
index $a$ such that $a_i-a_{i-1} > 1$ for at least one $i.$ Thus only the
indices $a$, which are Schubert symbols of certain Schubert cells of
$G_k(\C^{n+1})$, can provide nontrivial groups $\bar H_*(B_a,\pm \R).$
Moreover, for any Schubert symbol $a$ the set $B_a$ is diffeomorphic to the
complex vector space of dimension $\sum_{i=1}^{n} i \cdot (a_i - a_{i-1})$.
This implies statement B.
\medskip

{\sc Definition.} For any finite-dimensional topological space $X,$ its $k$-th
{\em self-join} $X^{*k}$ is defined as follows: we embed $X$ generically into
the space $\R^N$  of a very large dimension and take the union of all
$(k-1)$-dimensional simplices spanning all $k$-tuples of points in $X$. (The
genericity of the embedding means that these simplices have only obvious
intersections, namely, the intersection set of any two simplices is some their
common face.)
\medskip

{\sc Lemma 3.} {If $k\ge 2$, then $H_*((S^2)^{*k},\R)\equiv 0$ in all
dimensions.}
\medskip

{\em Proof.} The space $(S^2)^{*k}$ is obviously filtered by similar
self-joins:
\begin{equation}
S^1 \equiv (S^2)^{*1} \subset (S^2)^{*2} \subset \cdots \subset (S^2)^{*(k-1)}
\subset (S^2)^{*k}. \label{sjoins}
\end{equation}
Consider the spectral sequence generated by this filtration. Any space $\phi_i
\sm \phi_{i-1} \equiv (S^2)^{*i} \sm (S^2)^{*(i-1)}$ of this filtration is
fibered over $B(S^2,i)$ with the fiber $\stackrel{\circ}{\Delta}^{i-1}$. This
bundle is nontrivial, its fibers change their orientations exactly over the
same loops, which act nontrivially on the local system $\pm \R$. Hence the
first term of this sequence is given by the formula
$${\mathcal E}^1_{p,q} \simeq
\bar H_{p+q}((S^2)^{*p} \sm (S^2)^{*(p-1)},\R) \simeq \bar
H_{q+1}(B(\CP^1,p),\pm \R).$$ By Lemma 2 there are exactly two such nontrivial
cells ${\mathcal E}^1_{p,q}$, namely ${\mathcal E}^1_{2,1} \simeq {\mathcal
E}^1_{1,1} \simeq \R$. The unique differential $d^1$ of our spectral sequence,
which can be nontrivial, acts between these cells; it is exactly the boundary
homomorphism $H_3((S^2)^{*2},S^2;\R) \to H_2(S^2,\R);$ It follows immediately
from the shape of the generator of the former group (i.e., from that of the
2-dimensional Schubert cell $B_a,$ $a=(1,2)$ in $B(S^2,2)$) that this
differential an isomorphism. This implies our Lemma.
\medskip

Denote by $\tilde B(\CP^2,k)$ the subset in $B(\CP^2,k),$ consisting of {\em
generic} configurations, i.e. of such that none three points do lie in the same
line.
\medskip

{\sc Lemma 4.} {\em There are isomorphisms
\begin{equation}
\label{b3} \bar H_*(B(\CP^2,3),\pm \R) \stackrel{\sim}{\longrightarrow} \bar
H_*(\tilde B(\CP^2,3),\pm \R),
\end{equation}
\begin{equation}
\label{b4} \bar H_*(B(\CP^2,4),\pm \R) \stackrel{\sim}{\longrightarrow} \bar
H_*(\tilde B(\CP^2,4),\pm \R),
\end{equation}
induced by identical embeddings. Namely, both groups (\ref{b4}) are trivial in
all dimensions, and groups (\ref{b3}) are trivial in all dimensions other than
6 and are isomorphic to $\R$ in dimension 6.}
\medskip

{\em Proof.} The space $B(\CP^2,3) \sm \tilde B(\CP^2,3)$ consists of
3-configurations lying on the same line. This space is obviously fibered over
$\CP^{2\vee}$ with fiber equal to $B(\CP^1,3).$ Applying Lemma 2 B) to this
fiber, we get the assertion about the groups (\ref{b3}). Similarly, the space
$B(\CP^2,4) \sm \tilde B(\CP^2,4)$ is the union of two pieces, consisting of
4-configurations, some three  (respectively all four) points of which lie in
the same line. Both these pieces are fibered over $\CP^{2\vee}$  with fibers
$B(\CP^1,3) \times \C^2$ (respectively, $B(\CP^1,4)$). By Lemma 2 B) both these
pieces are $\pm \R$-acyclic, and Lemma 4 is proved.

\subsection{Proof of Theorem 2.}

There are exactly the following possible singular sets $K\subset \CP^2$ defined
by homogeneous polynomials $\C^3 \to \C^1$ of degree 3:

A) any point $u \in \CP^2;$

B) any pair of points in $\CP^2;$ the corresponding space $L(K)$ consists of
all 3-forms splitting into the product of two polynomials of degree 1 and 2,
vanishing at these two points;

C) any line in $\CP^2:$ the corresponding space $L(K)$ consists of 3-forms
vanishing on this line with multiplicity $\ge 2$;

D) any triple of points in $\CP^2$ not lying on the same line, i.e. defining a
point of the space $\tilde B(\CP^2,3)$. The space $L(K)$ in this case is
1-dimensional and consists of forms vanishing on any of three lines containing
any two of our three points;

E) the entire $\CP^2.$
\medskip

The corresponding terms $F_i \setminus F_{i-1}$ of our filtration and the
corresponding terms of the spectral sequence look as follows. Any of them is
the space of a fibered product of two fiber bundles, whose base is the
corresponding space of all possible sets $K$ (i.e., respectively, $\CP^2$,
$B(\CP^2,2)$, $\CP^{2\vee},$ $\tilde B(\CP^2,3)$ and $\{pt\}$). The fiber of
the first (respectively, second) bundle over any point $\{K\}$ of such a base
is the space $L(K)$ (respectively, the open cone
$\stackrel{\circ}{\Lambda}_{3,2}(K)$). The spaces of only second fiber bundles
are exactly the terms $\Phi_i \sm \Phi_{i-1}$ of the related filtration of the
continuous order complex $\Lambda_{3,2}$ of all possible sets $K$ defined by
homogeneous polynomials of order 3.

Let us describe explicitly all these terms and bundles.

A) The space $F_1 \equiv \tilde F_3$ is the space of a fiber bundle with fiber
$\C^7$ over $\CP^2$, this implies statement of Theorem 2 concerning the column
$p=1$\footnote{the fibers of the second bundle in this case are the points}.
\medskip

B) The space $F_2 \sm F_1 \equiv \tilde F_6 \sm \tilde F_5$. The fiber of the
second bundle over a point $(u,v) \in B(\CP^2,2)$ is the union of two segments,
joining this point $(u,v)$ with both points $u$ and $v\in \CP^2$, and not
containing their boundary points $u$ and $v$ (which belong to $F_1$). Obviously
such an union is homeomorphic to the open segment $(-1,1)$. Our bundle of such
segments is not oriented: it changes the orientation simultaneously with the
sheaf $\pm \Z$. Therefore the term $E^1_{2,q}$ of our spectral sequence is
given by the formula
\begin{equation}
\label{ss3_2} E^1_{2,q} \cong \bar H_{2+q}(F_2 \sm F_1, \R) \simeq \bar
H_{2+q-8}(\Phi_2 \sm \Phi_1,\R) \simeq \bar H_{2+q-9}(B(\CP^2,2),\pm \R).
\end{equation}

Thus the assertion of Theorem 2 concerning the column $E^1_{2,q}$ of the main
spectral sequence follows from Lemma 2.
\medskip

C) Let $K=\CP^1 \subset \CP^2.$ Then the reduced order complex
$\Lambda_{3,2}(K)$ is homeomorphic to the cone with vertex $\{L(K)\} \subset
\Omega_7,$ and base lying in $\Phi_2.$ By the  construction, this base
$\partial \Lambda_{3,2}(K)$ is canonically homeomorphic to the self-join
$(\CP^1)^{*2},$  whose homology groups are trivial by Lemma 3.

The space $F_3 \sm F_2$  is fibered over $\CP^{2\vee}$; its fiber over a line
$\{K\} \in \CP^{2\vee}$ is the product of the space $\C^3$ and the space
$\Lambda_{3,2}(K) \sm \Phi_2$, i.e. an open cone with base $\partial
\Lambda_{3,2}(K)$. Thus the Borel--Moore homology group of this fiber is
acyclic over $\R,$ this implies assertion of Theorem 2 concerning the column
$p=3.$
\medskip

D) $\Phi_4 \sm \Phi_3$  is fibered over $\tilde B(\CP^2,3);$ its fiber over a
point $(u,v,w) \subset \CP^2$ is the open triangle $\dot \Delta^2,$ whose
vertices are identified with these points $u,v$ and $w.$ In particular,
\begin{equation}
\label{d} \bar H_i(F_4 \sm F_3) \equiv \bar H_{i-2}(\Phi_4 \sm \Phi_3) \simeq
\bar H_{i-4}(\tilde B(\CP^2,3),\pm \Z).
\end{equation}

By Lemmas 4 and 2 B), the group $H_{j}(\tilde B(\CP^2,3),\pm \R)$ is isomorphic
to $\R$ if $j=6$ and is trivial for all  $j \ne 6$. This proves the assertion
of Theorem 2 concerning the column $p=4$.
\medskip

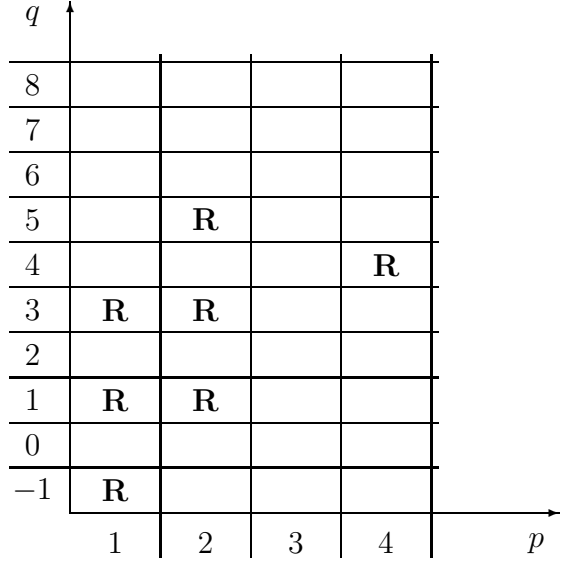
\begin{figure}
\unitlength=1.00mm \thinlines \linethickness{0.4pt}
\begin{picture}(75.00,78.00)
\put(10.00,10.00){\vector(1,0){65.00}} \put(72.00,6.00){\makebox(0,0)[cc]{$p$}}
\put(10.00,10.00){\vector(0,1){68.00}} \put(5.00,76.00){\makebox(0,0)[cc]{$q$}}
\put(5.00,67.00){\makebox(0,0)[cc]{$8$}}
\put(5.00,61.00){\makebox(0,0)[cc]{$7$}}
\put(5.00,55.00){\makebox(0,0)[cc]{$6$}}
\put(5.00,49.00){\makebox(0,0)[cc]{$5$}}
\put(5.00,43.00){\makebox(0,0)[cc]{$4$}}
\put(5.00,37.00){\makebox(0,0)[cc]{$3$}}
\put(5.00,31.00){\makebox(0,0)[cc]{$2$}}
\put(5.00,25.00){\makebox(0,0)[cc]{$1$}}
\put(5.00,19.00){\makebox(0,0)[cc]{$0$}}
\put(5.00,13.00){\makebox(0,0)[cc]{$-1$}} \put(22.00,4.00){\line(0,1){67.00}}
\put(34.00,71.00){\line(0,-1){67.00}} \put(46.00,4.00){\line(0,1){67.00}}
\put(58.00,71.00){\line(0,-1){67.00}} \put(16.00,37.00){\makebox(0,0)[cc]{${\bf
R}$}} \put(16.00,25.00){\makebox(0,0)[cc]{${\bf R}$}}
\put(16.00,13.00){\makebox(0,0)[cc]{${\bf R}$}}
\put(28.00,49.00){\makebox(0,0)[cc]{${\bf R}$}}
\put(28.00,37.00){\makebox(0,0)[cc]{${\bf R}$}}
\put(28.00,25.00){\makebox(0,0)[cc]{${\bf R}$}}
\put(52.00,43.00){\makebox(0,0)[cc]{${\bf R}$}}
\put(52.00,6.00){\makebox(0,0)[cc]{$4$}}
\put(40.00,6.00){\makebox(0,0)[cc]{$3$}}
\put(28.00,6.00){\makebox(0,0)[cc]{$2$}}
\put(16.00,6.00){\makebox(0,0)[cc]{$1$}} \put(2.00,16.00){\line(1,0){57.00}}
\put(2.00,22.00){\line(1,0){57.00}} \put(2.00,28.00){\line(1,0){57.00}}
\put(2.00,34.00){\line(1,0){57.00}} \put(2.00,40.00){\line(1,0){57.00}}
\put(2.00,46.00){\line(1,0){57.00}} \put(2.00,52.00){\line(1,0){57.00}}
\put(2.00,58.00){\line(1,0){57.00}} \put(2.00,64.00){\line(1,0){57.00}}
\put(2.00,70.00){\line(1,0){57.00}}
\end{picture}
\caption{spectral sequence for the term $F_5 \setminus F_4$.} \label{sf5}
\end{figure}

E) Finally, the set $F_5 \sm F_4$ coincides with $\Phi_5 \sm \Phi_4 \equiv
\Lambda_{3,2} \sm \Phi_4$; it is an open cone over the reduced topological
order complex $\Phi_4,$ so that $\bar H_i(F_5 \sm F_4) \simeq \tilde
H_{i-1}(\Phi_4)$. Consider the spectral sequence $e^a_{b,c} \to
H_*(\Phi_4,\R),$ generated by the filtration $\{\Phi_i\},$ $i=1,2,3,4.$ Its
term $e^1_{i,c}$ is isomorphic to $\bar H_{i+c}(\Phi_i \sm \Phi_{i-1},\R).$ The
homology groups of spaces $(\Phi_i \sm \Phi_{i-1})$ are related by the Thom
isomorphism with these of corresponding spaces $F_i \sm F_{i-1}$ (calculated in
the previous steps A --- D), hence the term $e^1$ of our sequence looks as is
shown in fig.~\ref{sf5}.

\medskip
{\sc Lemma 5.} {\em The differentials $\partial^1: e^1_{2,1} \to e^1_{1,1}$ and
$\partial^1: e^1_{2,3} \to e^1_{1,3}$ of our spectral sequence are
isomorphisms, as well as the map $\partial^2: e^2_{4,4} \to e^2_{2,5}$. In
particular, the term $e^\infty \equiv e^3$ of our spectral sequence consists of
unique nontrivial cell $e^3_{1,-1} \simeq \R$.}
\medskip

{\em Proof.} If some of these maps is not isomorphic, then there exists a
nontrivial group $H_{i-1}(\Phi_4,\R) \simeq \bar H_{i}(F_5 \sm F_4,\R)$ with
$i<9$. It follows from the shape of terms $E^1_{p,q}$ of the main spectral
sequence with $p \le 4$ (determined in the previous steps) that the
corresponding group $E^1_{5,i-5}$ cannot disappear in the further steps of this
sequence and thus the group $\bar H_i(\Sigma_{3,2}) \simeq H^{19-i}(\Pi_{3,2}
\sm \Sigma_{3,2})$ also is nontrivial. This is impossible because $\Pi_{3,2}
\sm \Sigma_{3,2}$ is a 10-dimensional Stein manifold. \quad $\Box$
\medskip

Thus we have proved the assertion of Theorem 2 concerning the shape of the term
$E^1$ of the main spectral sequence.
\medskip

{\sc Lemma 6.} {\em The homomorphism $d^1: E^1_{2,13} \to E^1_{1,13}$ is
trivial.}
\medskip

{\em Proof.} It is clear from fig.~\ref{ss3}, that the group $E_{1,15}$ cannot
disappear, hence $\bar H_{16}(\Sigma_{3,2},\R) \simeq H^3(\Pi_{3,2} \sm
\Sigma_{3,2},\R) \sim \R.$ $\Pi_{3,2} \sm \Sigma_{3,2}$ is a direct product of
something and $\C^*$, therefore the Poincar\'e polynomial of $ H^*(\Pi_{3,2}
\sm \Sigma_{3,2},\R)$ is divisible by $(1+t).$ Since $\bar
H_{17}(\Sigma_{3,2},\R) \simeq H^2(\Pi_{3,2} \sm \Sigma_{3,2},\R)=0,$ this
implies that $\bar H_{15}(\Sigma_{3,2},\R) \ne 0.$ But the unique cell
$E^1_{p,q}$ with $p+q=15$ is $E^1_{2,13},$ and our lemma is proved.
\medskip

It terminates the proof of Theorems 2 and 1.

\section{Non-singular quartics in $\CP^2$}

{\sc Theorem 3.} {\em The Poincar\'e polynomial of the real cohomology group of
the space $N_{4,2}$ of nonsingular curves of degree 4 in $\CP^2$ is equal to
$(1+t^3)(1+t^5)(1+t^6)$.}
\medskip

For the proof, consider the standard conical resolution of $\Sigma_{4,2}$
described in \S~3. To describe it explicitly, we must classify all singular
sets in $\CP^2$ of polynomials $\C^3 \to \C^1$ of degree 4.
\medskip

{\sc Proposition 6.} {\em There are exactly the following possible singular
sets in $\CP^2$ of homogeneous polynomials of degree 4 in $\C^3$ (in angular
parentheses we indicate the complex dimension of the space of all polynomials
having the singularity at an arbitrary fixed set of the corresponding class):

\begin{enumerate}
\item
Any point in $\CP^2$ \quad $\langle 12 \rangle$

\item
Any pair of points \quad $\langle 9 \rangle$

\item
Any three point on the same line \quad $\langle 7 \rangle$ \unitlength=1.00mm
\thinlines \linethickness{0.4pt}
\begin{picture}(17.67,6.00)
\put(7.00,3.00){\oval(6.00,6.00)[b]} \put(13.00,3.00){\oval(6.00,6.00)[t]}
\put(17.67,5.67){\line(-5,-2){15.33}}
\end{picture}

\item
Any generic triple of points (i.e. not lying on the same line) \quad $\langle 6
\rangle$

\item
Any line $\CP^1 \subset \CP^2$ \quad $\langle 6 \rangle$

\item
Any three points on the same line plus a point not on the line \quad $\langle 4
\rangle$ \unitlength=1.00mm \thinlines \linethickness{0.4pt}
\begin{picture}(11.00,9.00)
\put(5.50,2.50){\oval(7.00,3.00)[r]} \put(5.33,5.00){\oval(6.00,8.00)[lb]}
\put(5.33,0.00){\oval(6.00,8.00)[lt]} \put(1.00,0.67){\line(5,2){10.00}}
\end{picture}

\item
Any generic quadruple of points (i.e., none three of them are on the same line;
the corresponding polynomial necessarily splits into the product of two
quadrics) \quad $\langle 3 \rangle$

\item
Any line $\CP^1 \subset \CP^2$ plus a point outside it \quad $\langle 3
\rangle$

\item
Any five points, four of which are generic, and the fifth is the intersection
point of two lines spanning some pairs of first four points (the corresponding
polynomial splits into the product of two linear functions and one quadric)
\quad $\langle 2 \rangle$ \unitlength=1.00mm \thinlines \linethickness{0.4pt}
\begin{picture}(9.00,5.00)
\put(6.00,2.00){\circle{6.00}} \put(3.00,-1.00){\line(1,1){6.00}}
\put(9.00,-1.00){\line(-1,1){6.00}}
\end{picture}

\item
Six points, which are all the intersection points of some four lines in general
position \quad $\langle 1 \rangle$ \unitlength=1.00mm \thinlines
\linethickness{0.4pt}
\begin{picture}(12.00,5.50)
\put(4.00,1.00){\line(1,0){8.00}} \put(6.00,-2.00){\line(0,1){7.50}}
\put(5.00,5.00){\line(1,-1){6.00}} \put(11.00,5.00){\line(-1,-1){7.00}}
\end{picture}

\item
Any nonsingular quadric in $\CP^2$ \quad $\langle 1 \rangle$

\item
Two lines in $\CP^2$ \quad $\langle 1 \rangle$

\item
Entire $\CP^2$. \quad $\langle 0 \rangle$ \quad $\Box$
\end{enumerate}
}

I have  liberty to improve slightly the canonical filtration in the resolution
space $\sigma_{4,2}$ of the discriminant and in the underlying order complex
$\Lambda_{4,2}$  of singular sets. Namely, for any $i=1, \ldots, 13,$ I denote
by $F_i$ (respectively, by $\Phi_i$) the union of all sets $\sigma_{4,2}(K)
\subset \sigma_{4,2}$ (respectively, $\Lambda_{4,2}(K) \subset \Lambda_{4,2}$)
where $K$ is one of sets described in items $1, \ldots, i$ of the previous
proposition.  \medskip

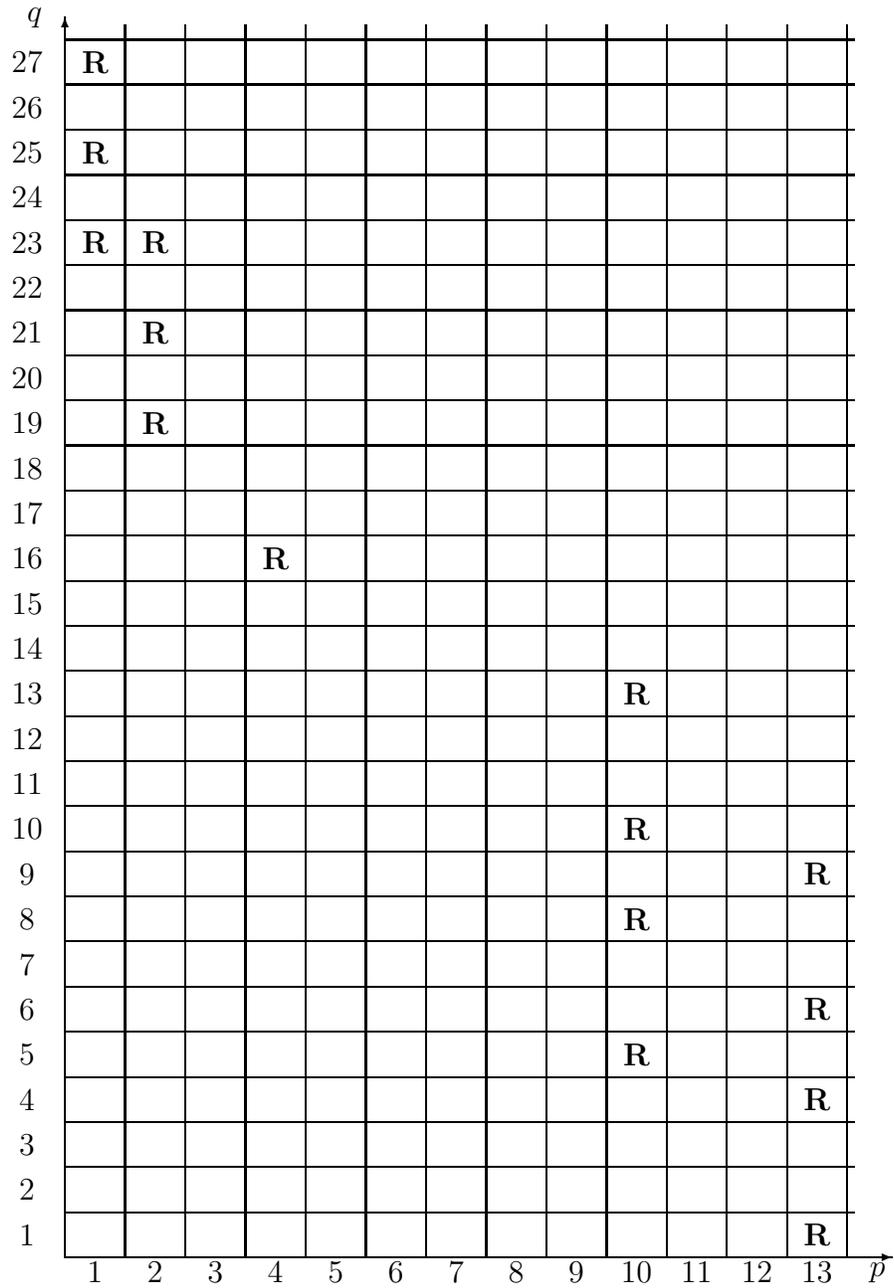
\begin{figure}
\unitlength=1.00mm \thinlines \linethickness{0.4pt}
\begin{picture}(120.00,162.00)(0,-3)
\put(10.00,-3.00){\vector(1,0){110.00}} \put(10.00,-3.00){\vector(0,1){165.00}}
\put(10.00,3.00){\line(1,0){105.00}} \put(10.00,9.00){\line(1,0){105.00}}
\put(10.00,15.00){\line(1,0){105.00}} \put(10.00,21.00){\line(1,0){105.00}}
\put(10.00,27.00){\line(1,0){105.00}} \put(10.00,33.00){\line(1,0){105.00}}
\put(10.00,39.00){\line(1,0){105.00}} \put(10.00,45.00){\line(1,0){105.00}}
\put(10.00,51.00){\line(1,0){105.00}} \put(10.00,57.00){\line(1,0){105.00}}
\put(10.00,63.00){\line(1,0){105.00}} \put(10.00,69.00){\line(1,0){105.00}}
\put(10.00,75.00){\line(1,0){105.00}} \put(10.00,81.00){\line(1,0){105.00}}
\put(10.00,87.00){\line(1,0){105.00}} \put(18.00,-3.00){\line(0,1){164.00}}
\put(26.00,-3.00){\line(0,1){164.00}} \put(34.00,-3.00){\line(0,1){164.00}}
\put(42.00,-3.00){\line(0,1){164.00}} \put(50.00,-3.00){\line(0,1){164.00}}
\put(58.00,-3.00){\line(0,1){164.00}} \put(66.00,-3.00){\line(0,1){164.00}}
\put(74.00,-3.00){\line(0,1){164.00}} \put(82.00,-3.00){\line(0,1){164.00}}
\put(90.00,-3.00){\line(0,1){164.00}} \put(98.00,-3.00){\line(0,1){164.00}}
\put(106.00,-3.00){\line(0,1){164.00}} \put(114.00,-3.00){\line(0,1){164.00}}
\put(14.00,-5.00){\makebox(0,0)[cc]{1}} \put(22.00,-5.00){\makebox(0,0)[cc]{2}}
\put(30.00,-5.00){\makebox(0,0)[cc]{3}} \put(38.00,-5.00){\makebox(0,0)[cc]{4}}
\put(46.00,-5.00){\makebox(0,0)[cc]{5}} \put(54.00,-5.00){\makebox(0,0)[cc]{6}}
\put(62.00,-5.00){\makebox(0,0)[cc]{7}} \put(70.00,-5.00){\makebox(0,0)[cc]{8}}
\put(78.00,-5.00){\makebox(0,0)[cc]{9}}
\put(86.00,-5.00){\makebox(0,0)[cc]{10}}
\put(94.00,-5.00){\makebox(0,0)[cc]{11}}
\put(102.00,-5.00){\makebox(0,0)[cc]{12}}
\put(110.00,-5.00){\makebox(0,0)[cc]{13}}
\put(5.00,90.00){\makebox(0,0)[cc]{16}} \put(5.00,96.00){\makebox(0,0)[cc]{17}}
\put(5.00,102.00){\makebox(0,0)[cc]{18}}
\put(5.00,108.00){\makebox(0,0)[cc]{19}}
\put(5.00,114.00){\makebox(0,0)[cc]{20}}
\put(5.00,120.00){\makebox(0,0)[cc]{21}}
\put(5.00,126.00){\makebox(0,0)[cc]{22}}
\put(5.00,132.00){\makebox(0,0)[cc]{23}}
\put(5.00,138.00){\makebox(0,0)[cc]{24}}
\put(5.00,144.00){\makebox(0,0)[cc]{25}}
\put(5.00,150.00){\makebox(0,0)[cc]{26}}
\put(5.00,156.00){\makebox(0,0)[cc]{27}}
\put(6.00,162.00){\makebox(0,0)[cc]{$q$}}
\put(118.00,-5.00){\makebox(0,0)[cc]{$p$}}
\put(10.00,93.00){\line(1,0){105.00}} \put(10.00,99.00){\line(1,0){105.00}}
\put(10.00,105.00){\line(1,0){105.00}} \put(10.00,111.00){\line(1,0){105.00}}
\put(10.00,117.00){\line(1,0){105.00}} \put(10.00,123.00){\line(1,0){105.00}}
\put(10.00,129.00){\line(1,0){105.00}} \put(10.00,135.00){\line(1,0){105.00}}
\put(10.00,141.00){\line(1,0){105.00}} \put(10.00,147.00){\line(1,0){105.00}}
\put(10.00,153.00){\line(1,0){105.00}} \put(10.00,159.00){\line(1,0){105.00}}
\put(14.00,156.00){\makebox(0,0)[cc]{${\bf R}$}}
\put(14.00,144.00){\makebox(0,0)[cc]{${\bf R}$}}
\put(14.00,132.00){\makebox(0,0)[cc]{${\bf R}$}}
\put(22.00,132.00){\makebox(0,0)[cc]{${\bf R}$}}
\put(22.00,120.00){\makebox(0,0)[cc]{${\bf R}$}}
\put(22.00,108.00){\makebox(0,0)[cc]{${\bf R}$}}
\put(38.00,90.00){\makebox(0,0)[cc]{${\bf R}$}}
\put(5.00,84.00){\makebox(0,0)[cc]{15}} \put(5.00,78.00){\makebox(0,0)[cc]{14}}
\put(5.00,72.00){\makebox(0,0)[cc]{13}} \put(5.00,66.00){\makebox(0,0)[cc]{12}}
\put(5.00,60.00){\makebox(0,0)[cc]{11}} \put(5.00,54.00){\makebox(0,0)[cc]{10}}
\put(5.00,48.00){\makebox(0,0)[cc]{9}} \put(5.00,42.00){\makebox(0,0)[cc]{8}}
\put(5.00,36.00){\makebox(0,0)[cc]{7}} \put(5.00,30.00){\makebox(0,0)[cc]{6}}
\put(5.00,24.00){\makebox(0,0)[cc]{5}} \put(5.00,18.00){\makebox(0,0)[cc]{4}}
\put(5.00,12.00){\makebox(0,0)[cc]{3}} \put(5.00,6.00){\makebox(0,0)[cc]{2}}
\put(5.00,0.00){\makebox(0,0)[cc]{1}} \put(10.00,3.00){\line(1,0){105.00}}
\put(10.00,9.00){\line(1,0){105.00}} \put(10.00,15.00){\line(1,0){105.00}}
\put(10.00,21.00){\line(1,0){105.00}} \put(86.00,72.00){\makebox(0,0)[cc]{${\bf
R}$}} \put(86.00,54.00){\makebox(0,0)[cc]{${\bf R}$}}
\put(86.00,42.00){\makebox(0,0)[cc]{${\bf R}$}}
\put(86.00,24.00){\makebox(0,0)[cc]{${\bf R}$}}
\put(110.00,48.00){\makebox(0,0)[cc]{${\bf R}$}}
\put(110.00,30.00){\makebox(0,0)[cc]{${\bf R}$}}
\put(110.00,18.00){\makebox(0,0)[cc]{${\bf R}$}}
\put(110.00,0.00){\makebox(0,0)[cc]{${\bf R}$}}
\end{picture}
\caption{Spectral sequence for the space of nonsingular homogeneous polynomials
of degree 4 in $\C^3$} \label{ssquart}
\end{figure}

{\sc Proposition 7.} {\em The term $E^1$ of the spectral sequence, generated by
this filtration and converging to the group $\bar H_*(\sigma_{4,2}, \R),$ is as
shown in Fig.~\ref{ssquart}. This spectral sequence degenerates at this term:
$E^\infty \equiv E^1$.}
\medskip

{\em Proof.} The shape of columns $p={\bf 1,2}$ and ${\bf 4}$ of $E^1$ is
justified in essentially the same way as that of the same columns of the
spectral sequence presented in Fig.~\ref{ss3}. The triviality of columns
$p={\bf 3}$ and $p={\bf 6}$ follows from Lemma 2 ($n=1,$ $k=3$), and that of
the column $p={\bf 7}$ from the same lemma and additionally Lemma 4. For any
set $K$ of type ${\bf 5}$  (i.e., $K = \CP^1$) the link $\partial
(\Lambda_{4,2}(K))$  of the corresponding order complex is the space
$(\CP^1)^{*3}$, thus the triviality of the column $p=5$ follows from Lemma 3.
Similarly, for any set $K$ of type ${\bf 11}$  (i.e., $K$ is a generic quadric
in $\CP^2$) the link $\partial (\Lambda_{4,2}(K))$ is the space $K^{*4}$, and
the triviality of the column $p=11$ also follows from Lemma 3.

The space of configurations of type ${\bf 9}$ is fibered over
$B(\CP^{2\vee},2)$ with fiber $B(\C^1,2) \times B(\C^1,2)$, and everything dyes
already over this fiber. Thus the column $p=9$ also is trivial.

In the case mentioned in the column $p={\bf 10}$ the term $F_{10} \sm F_9$ is
the space of a fiber bundle, whose base is the configuration space $\tilde
B(\CP^{2\vee},4)$ of all generic collections of four lines in $\CP^2$ (i.e.,
such that none three of them intersect at the same point), and the fiber is the
direct product of a complex line and a $5$-dimensional open simplex, whose 6
vertices correspond to all intersection points of these four lines. It is easy
to check that this bundle of 5-simplices is {\em orientable}, thus $\bar
H_*(F_{10}\sm F_9,\R) \simeq \bar H_{*-7}(\tilde B(\CP^{2\vee},4),\R).$ Let us
calculate the latter homology group. By the projective duality, $\tilde
B(\CP^{2\vee},4)$ is diffeomorphic to the space $\tilde B(\CP^{2},4)$ of all
generic collections of 4 points, see e.g. Lemma 4. The similar space $\tilde
{\mathcal F}(\CP^2,4)$  of {\em ordered} generic collections of points is
obviously diffeomorphic to the group $PGL(\CP^2),$ whose real homology groups
are well-known. It is easy to check that these homology groups do not disappear
when we pass to the base of the 24-fold covering $\tilde {\mathcal F}(\CP^2,4)
\stackrel{S(4)}{\longrightarrow} \tilde B(\CP^2,4)$, and the shape of the
column $p=10$ follows.

The triviality of the column $p={\bf 8}$ follows from the following lemma.
\medskip

{\sc Lemma 7.} {\em For any set $K \subset \CP^2$ of type 8, $K =\CP^1 \sqcup
*$, the link $\partial \Lambda_{4,2}(K)$ of the corresponding order complex
$\Lambda_{4,2}(K)$ is homeomorphic to the suspension $\Sigma(\partial
\Lambda_{4,2}(\CP^1))$ of the similar link for the set $\CP^1$ of type 5.}
\medskip

{\em Proof.} This link consists of two pieces, the first of which is the order
complex $\Lambda_{4,2}(\CP^1)$ (obviously homeomorphic to the cone over its own
link), and the second is the union of all order complexes
$\Lambda_{4,2}(\kappa),$ where $\kappa$  is some subset of type 6 in $K$. All
the latter order complexes are the tetrahedra, one of whose 4 vertices is the
distinguished point $* = K \sm \CP^1,$ and another three are arbitrary
different points of $\CP^1$. Thus the second piece also is homeomorphic to the
cone over $\partial \Lambda_{4,2} \sim (S^2)^{*3}$. The intersection of these
two cones is their common base, and Lemma 7 is proved.
\medskip

Now let $K$ be of type ${\bf 12}$,  i.e. the union of two complex lines $l_1,
l_2 \subset \CP^2.$ The corresponding link $\partial \Lambda_{4,2}(K)$ is
covered by two diffeomorphic sets, $A_1$ and $A_2$. Namely, $A_i$, $i=1,2,$ is
the union of all continuous order complexes $\Lambda_{4,2}(\kappa),$ where
$\kappa$ is the set of type 8 consisting of the line $l_i$ and a point of the
other line $K \sm l_i.$ The triviality of the column $p=12$ follows now from
the Mayer--Vietoris formula and the following lemma.
\medskip

{\sc Lemma 8.} {\em All three complexes $A_1$, $A_2$ and $A_1 \cap A_2$ are
acyclic (over $\R$) in all positive dimensions.}
\medskip

{\em Proof.} Consider the following filtration of the set $A_1$:
$$\Lambda_{4,2}(l_1) \subset A_1 \cap \Phi_6 \subset A_1 \cap \Phi_8
\equiv A_1.$$ The set $A_1 \sm (A_1 \cap \Phi_6)$ is a fiber bundle, whose base
is the space $(K \sm l_1) \sim \C^1$, parametrizing all possible sets of type 8
containing $l_1,$ and fiber $\stackrel{\circ}{\Lambda}_{4,2}(\kappa),$ where
$\kappa$ is such a set. By Lemma 7 the Borel--Moore homology group of this
difference (or, which is the same, the relative group $H_*(A_1, (A_1 \cap
\Phi_6);\R)$) is trivial.

The space $(A_1 \cap \Phi_6) \sm \Lambda_{4,2}(l_1)$ also is a fiber bundle
with base $(K \sm l_1) \sim \C^1$; its fiber is the space $\partial
\Lambda_{4,2}(l_1) \sim (S^2)^{*3}$. By Lemma 3 the Borel--Moore homology group
of this fiber also is trivial. Thus $H_*(A_1) $ coincides with the homology
group of the space $ \Lambda_{4,2}(l_1),$ which is contractible, and acyclicity
of $A_1$ is completely proved. The proof for $A_2$ is exactly the same.

Finally, $A_1 \cap A_2$ is the continuous order complex of following sets:

a) any point of $K$;

b) any pair of points in $K$ lying in different lines (one of which can
coincide with the intersection point $l_1 \cap l_2$);

c) any triple of points in $K$, the first of which is the crossing point $l_1
\cap l_2$, and two other lie in different lines $l_1,$ $l_2.$

The homology groups of this complex can be easily calculated and coincide with
the homology groups of a point. Lemma 8 is completely proved and the column
$p=12$ actually is empty. \medskip

As in the proof of Theorem 2, the ultimate column $p=13$ is calculated by a
spectral sequence, whose twelve columns repeat these of the main sequence, with
coordinates $q$ decreased by twice the number shown in angular parentheses in
Proposition 6. The fact that its lower terms kill one another is proved in
exactly the same way as it was done in the proof of Theorem 2, see Lemma 5.
(The complement of the discriminant is a Stein manifold and cannot have too
low-dimensional homology groups).

Finally, the fact that the term $E^1$, shown in Fig.~\ref{ssquart}, coincides
with the stable term $E^\infty,$ also is proved similarly to the similar
statement (Lemma 6) in the proof of Theorem 2. Say, the cell $E_{2,23}$  cannot
disappear, because it should match the cell $E_{1,26}$, which, in its turn,
cannot disappear by dimensional reasons. Cell $E_{13,9}$ cannot disappear,
because it should match one of cells $E_{2,21}$ or $E_{10,13}$, etc. \quad
$\Box$

\section{Non-singular cubical surfaces in $\CP^3$}

{\sc Theorem 4.} {\em The Poincar\'e polynomial of the real cohomology group of
the space $N_{3,3}$ of nonsingular cubical hypersurfaces in $\CP^3$ is equal to
$(1+t^3)(1+t^5)(1+t^7)$. In particular, again $H^*(N_{3,3},\R) \simeq
H^*(PGL(\CP^3),\R).$}
\medskip

{\sc Proposition 8.} {\em There are exactly the following possible singular
sets in $\CP^3$ of homogeneous cubical polynomials in $\C^4$:

\begin{enumerate}
\item
Any point in $\CP^3$ \quad $\langle 16 \rangle$

\item
Any pair of points \quad $\langle 12 \rangle$

\item
Any line \quad $\langle 10 \rangle$

\item
Any triple of points in general position (i.e. not on the same line) \quad
$\langle 8 \rangle$

\item \label{plquad}
Any generic quadric inside any plane $\CP^2 \subset \CP^3$ \quad $\langle 5
\rangle$

\item
Any pair of intersecting lines in $\CP^3$ \quad $\langle 5 \rangle$

\item
Any generic quadruple of points (i.e. not lying in the same two-dimensional
plane) \quad $\langle 4 \rangle$

\item Any plane $\CP^2 \subset \CP^3$
\quad $\langle 4 \rangle$

\item
Any triple of lines having one common point (such a polynomial is equal to $F=
xyz$ in appropriate coordinates) \quad $\langle 1 \rangle$

\item
Any generic plane quadric (as in item \ref{plquad}) plus any point not in the
same plane \quad $\langle 1 \rangle$

\item
Entire $\CP^3.$ \quad $\langle 0 \rangle$
\end{enumerate}
}

\begin{figure}
\unitlength=1.00mm \thinlines \linethickness{0.4pt}
\begin{picture}(130.00,145.00)
\put(10.00,5.00){\vector(0,1){140.00}} \put(10.00,5.00){\vector(1,0){113.00}}
\put(99.00,11.00){\line(-1,0){89.00}} \put(99.00,17.00){\line(-1,0){89.00}}
\put(99.00,23.00){\line(-1,0){89.00}} \put(99.00,29.00){\line(-1,0){89.00}}
\put(99.00,35.00){\line(-1,0){89.00}} \put(99.00,41.00){\line(-1,0){89.00}}
\put(99.00,47.00){\line(-1,0){89.00}} \put(99.00,53.00){\line(-1,0){89.00}}
\put(99.00,59.00){\line(-1,0){89.00}} \put(99.00,65.00){\line(-1,0){89.00}}
\put(99.00,71.00){\line(-1,0){89.00}} \put(99.00,77.00){\line(-1,0){89.00}}
\put(99.00,83.00){\line(-1,0){89.00}} \put(99.00,89.00){\line(-1,0){89.00}}
\put(99.00,95.00){\line(-1,0){89.00}} \put(99.00,101.00){\line(-1,0){89.00}}
\put(99.00,107.00){\line(-1,0){89.00}} \put(99.00,113.00){\line(-1,0){89.00}}
\put(99.00,119.00){\line(-1,0){89.00}} \put(99.00,125.00){\line(-1,0){89.00}}
\put(99.00,131.00){\line(-1,0){89.00}} \put(99.00,137.00){\line(-1,0){89.00}}
\put(5.00,8.00){\makebox(0,0)[cc]{16}} \put(5.00,14.00){\makebox(0,0)[cc]{17}}
\put(5.00,20.00){\makebox(0,0)[cc]{18}} \put(5.00,26.00){\makebox(0,0)[cc]{19}}
\put(5.00,32.00){\makebox(0,0)[cc]{20}} \put(5.00,38.00){\makebox(0,0)[cc]{21}}
\put(5.00,44.00){\makebox(0,0)[cc]{22}} \put(5.00,50.00){\makebox(0,0)[cc]{23}}
\put(5.00,56.00){\makebox(0,0)[cc]{24}} \put(5.00,62.00){\makebox(0,0)[cc]{25}}
\put(5.00,68.00){\makebox(0,0)[cc]{26}} \put(5.00,74.00){\makebox(0,0)[cc]{27}}
\put(5.00,80.00){\makebox(0,0)[cc]{28}} \put(5.00,86.00){\makebox(0,0)[cc]{29}}
\put(5.00,92.00){\makebox(0,0)[cc]{30}} \put(5.00,98.00){\makebox(0,0)[cc]{31}}
\put(5.00,104.00){\makebox(0,0)[cc]{32}}
\put(5.00,110.00){\makebox(0,0)[cc]{33}}
\put(5.00,116.00){\makebox(0,0)[cc]{34}}
\put(5.00,122.00){\makebox(0,0)[cc]{35}}
\put(5.00,128.00){\makebox(0,0)[cc]{36}}
\put(5.00,134.00){\makebox(0,0)[cc]{37}} \put(18.00,5.00){\line(0,1){133.00}}
\put(26.00,5.00){\line(0,1){133.00}} \put(34.00,5.00){\line(0,1){133.00}}
\put(42.00,5.00){\line(0,1){133.00}} \put(50.00,5.00){\line(0,1){133.00}}
\put(58.00,5.00){\line(0,1){133.00}} \put(66.00,5.00){\line(0,1){133.00}}
\put(74.00,5.00){\line(0,1){133.00}} \put(82.00,5.00){\line(0,1){133.00}}
\put(90.00,5.00){\line(0,1){133.00}} \put(98.00,5.00){\line(0,1){133.00}}
\put(14.00,2.00){\makebox(0,0)[cc]{1}} \put(22.00,2.00){\makebox(0,0)[cc]{2}}
\put(30.00,2.00){\makebox(0,0)[cc]{3}} \put(38.00,2.00){\makebox(0,0)[cc]{4}}
\put(46.00,2.00){\makebox(0,0)[cc]{5}} \put(54.00,2.00){\makebox(0,0)[cc]{6}}
\put(62.00,2.00){\makebox(0,0)[cc]{7}} \put(70.00,2.00){\makebox(0,0)[cc]{8}}
\put(78.00,2.00){\makebox(0,0)[cc]{9}} \put(86.00,2.00){\makebox(0,0)[cc]{10}}
\put(94.00,2.00){\makebox(0,0)[cc]{11}}
\put(14.00,134.00){\makebox(0,0)[cc]{${\bf R}$}}
\put(14.00,122.00){\makebox(0,0)[cc]{${\bf R}$}}
\put(14.00,110.00){\makebox(0,0)[cc]{${\bf R}$}}
\put(14.00,98.00){\makebox(0,0)[cc]{${\bf R}$}}
\put(38.00,68.00){\makebox(0,0)[cc]{${\bf R}$}}
\put(38.00,56.00){\makebox(0,0)[cc]{${\bf R}$}}
\put(38.00,44.00){\makebox(0,0)[cc]{${\bf R}$}}
\put(38.00,32.00){\makebox(0,0)[cc]{${\bf R}$}}
\put(22.00,110.00){\makebox(0,0)[cc]{${\bf R}$}}
\put(22.00,98.00){\makebox(0,0)[cc]{${\bf R}$}}
\put(22.00,74.00){\makebox(0,0)[cc]{${\bf R}$}}
\put(22.00,62.00){\makebox(0,0)[cc]{${\bf R}$}}
\put(22.00,86.00){\makebox(0,0)[cc]{${\bf R}^2$}}
\put(62.00,8.00){\makebox(0,0)[cc]{${\bf R}$}}
\put(115.00,2.00){\makebox(0,0)[cc]{$p$}}
\put(5,142.00){\makebox(0,0)[cc]{$q$}}
\end{picture}
\caption{Spectral sequence for nonsingular cubical polynomials in $\C^4$}
\label{sssurf}
\end{figure}
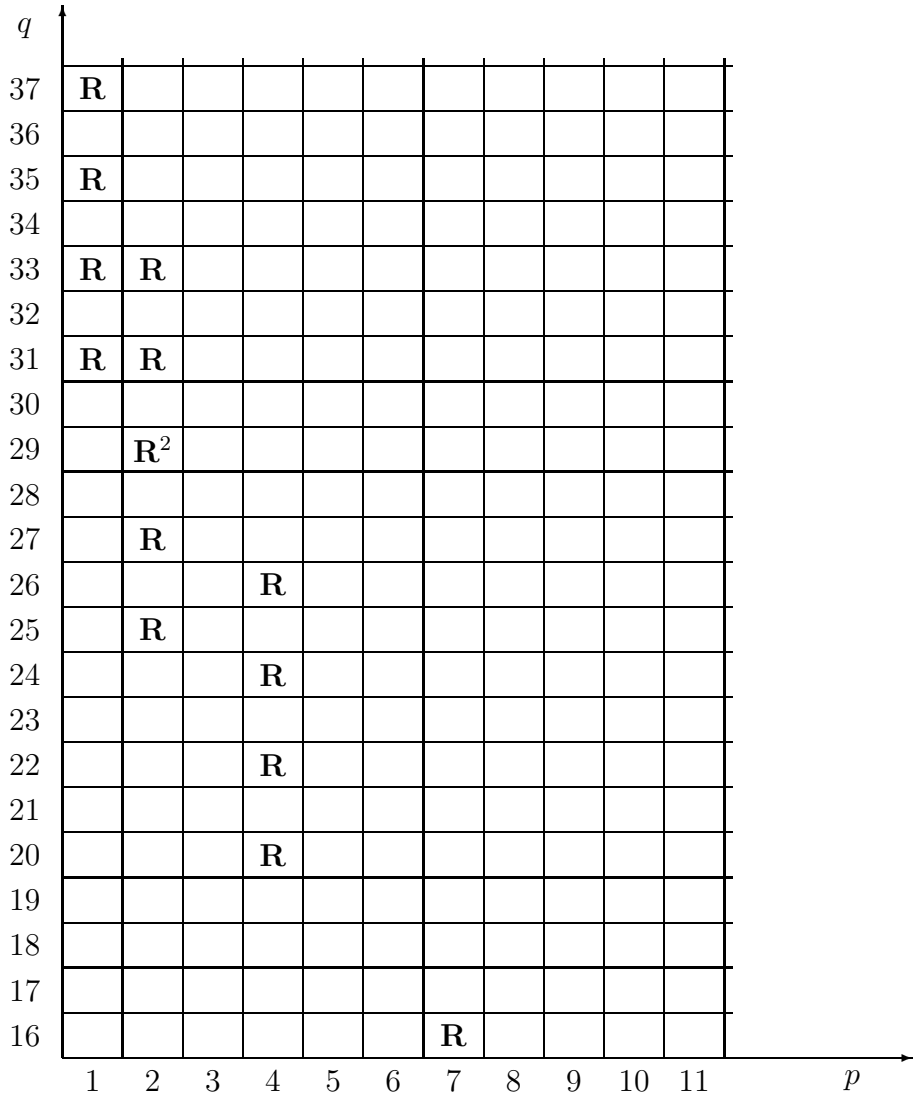

Exactly as in the proof of Theorem 3, for any $i=1, \ldots, 11,$ denote by
$F_i$ the union of all sets $\sigma_{3,3}(K) \subset \sigma_{3,3}$ where $K$ is
one of sets described in items $1, \ldots, i$ of the previous proposition.
\medskip

{\sc Proposition 9.} {\em The term $E^1$ of the spectral sequence, generated by
this filtration and converging to the group $\bar H_*(\sigma_{3,3}, \R),$ is as
shown in Fig.~\ref{sssurf}. This spectral sequence degenerates at this term:
$E^\infty \equiv E^1$.}
\medskip

{\em Proof.} The shape of columns $p={\bf 1}$ and ${\bf 2}$ follows immediately
from Lemma 2B. The justification of columns $p={\bf 4}$ and ${\bf 7}$ needs
additionally the following lemma.

Denote by $\tilde B(\CP^3,k)$ the subset in $B(\CP^3,k),$ consisting of {\em
generic} configurations, i.e. of such that none their three points lie on the
same line and none four lie on the same 2-plane.
\medskip

{\sc Lemma 9.} {\em There are isomorphisms
\begin{equation}
\label{b33} \bar H_*(B(\CP^3,3),\pm \R) \stackrel{\sim}{\longrightarrow} \bar
H_*(\tilde B(\CP^3,3),\pm \R),
\end{equation}
\begin{equation}
\label{b43} \bar H_*(B(\CP^3,4),\pm \R) \stackrel{\sim}{\longrightarrow} \bar
H_*(\tilde B(\CP^3,4),\pm \R),
\end{equation}
induced by identical embeddings. Namely, both groups (\ref{b33}) are isomorphic
to $\R$ in dimensions 6,8,10 and 12, and are trivial in all other dimensions,
and groups (\ref{b43}) are trivial in all dimensions other than 12 and are
isomorphic to $\R$ in dimension 12.}
\medskip

{\em Proof.} The isomorphism (\ref{b33}) follows from exactly the same
arguments as (\ref{b3}), only with space $\CP^{2\vee}$ of all lines in $\CP^2$
replaced by the space $G_2(\C^4)$ of lines in $\CP^3$.

The space $B(\CP^3,4) \sm \tilde B(\CP^3,4)$ is the union of two pieces,
consisting of 4-configurations, whose projective span is a 2-plane
(respectively, a line) in $\CP^3$.  By Lemmas 2 B) and 4 the $\pm \R$-homology
groups of both these pieces are trivial, and Lemma 9 is proved.
\medskip

For $K$ of type ${\bf 3}$, $K=\CP^1 \subset \CP^2,$ $\partial \Lambda_{3,3}(K)
\simeq (S^2)^{*2},$  therefore the shape of the column $p={\bf 3}$ follows from
Lemma 3.

Similarly, for $K$ of type ${\bf 5}$ $\partial \Lambda_{3,3}(K) \simeq
(S^2)^{*3},$ and the shape of the column $p=5$ follows from the same lemma.

For $K$ of type ${\bf 6}$,  $K=l_1 \cup l_2,$ the link $\partial
\Lambda_{3,3}(K)$ is covered by three subsets: $A_1=\Lambda_{3,3}(l_1),$
$A_2=\Lambda_{3,3}(l_2),$ and $A_3$ swept out by all order complexes
$\Lambda_{3,3}(\kappa),$ where $\kappa$ are 3-point sets, one whose point is
the crossing point $l_1 \cap l_2,$ and two other lie in different lines $l_1$
and $l_2$. The order complexes $A_1$ and $A_2$ are contractible, and their
intersection consists of one point, thus $H_*(A_1 \cup A_2) \simeq H_*(pt).$
Further, $A_3 \sm (A_1 \cup A_2)$  is a fiber bundle, whose base is the space
$\C^1 \times \C^1 \equiv (K \sm l_1) \times (K \sm l_2)$ parametrizing all
above-described sets $\kappa$, and the fiber over such a set is the triangle
with two sides removed. Thus the Borel--Moore homology group of such a fiber is
trivial, and by the K\"unneth formula
$$\bar H_*(\partial \Lambda_{3,3}(K) \sm (A_1 \cup A_2)) \equiv
H_*(\partial \Lambda_{3,3}(K), (A_1 \cup A_2)) \equiv 0.$$ This proves the
triviality of the column $p=6.$

Similarly, if $K$ is of type ${\bf 9},$ $K = l_1 \cap l_2 \cap l_3,$ then the
link $\partial \Lambda_{3,3}(K)$ is covered by four sets $A_1, A_2, A_3, A_4,$
where $A_i,$ $i=1,2,$and $3,$ are the order complexes $\Lambda_{3,3}(\mu_i),$
$\mu_i$ are unions of some two of three lines constituting $K$. $A_4$ is swept
out by
 all order complexes $\Lambda_{3,3}(\kappa),$
where $\kappa$ are 4-point sets, one whose point is the crucial point of $K$,
and three other lie in different lines $l_1, l_2$ and $l_3$. $A_4 \sm (A_1 \cup
A_2\cup A_3)$  is a fiber bundle, whose base is the space $\C^1 \times \C^1
\times \C^1$ parametrizing all above-described sets $\kappa$, and the fiber
over such a set is the tetrahedron with three maximal faces removed. Thus the
Borel--Moore homology group of such a fiber is trivial, and by the K\"unneth
formula
$$\bar H_*(\partial \Lambda_{3,3}(K) \sm (A_1 \cup A_2 \cup A_3)) \equiv
H_*(\partial \Lambda_{3,3}(K), (A_1 \cup A_2 \cup A_3)) \equiv 0.$$ Further,
the order complexes $A_1,$ $A_2$ and $A_3$ are acyclic, their pairwise
intersections are the order complexes $\Lambda_{3,3}(l_j),$ which also are
acyclic, and their total intersection is a point. This proves the triviality of
the column $p=9.$

For $K$ of type ${\bf 10},$ (i.e. consisting of a generic plane quadric $\tilde
K$ and one point not in the same plane), the link $\partial \Lambda_{3,3}(K)$
is homeomorphic to the suspension $\Sigma (\partial \Lambda_{3,3}(\tilde K)),$
cf. Lemma 7. This proves the triviality of the column $p=10.$

{\em Column} $p={\bf 8}$. Let $K = \CP^2 \subset \CP^3.$ The link $\partial
\Lambda_{3,3}(K)$ is filtered by its intersections with terms  $\Phi_i$ of the
filtration of $\Lambda_{3,3}.$ We need only the following segment of this
filtration:
\begin{equation}
\label{filtr}
\partial \Lambda_{3,3}(K) \cap \Phi_4 \subset
\partial \Lambda_{3,3}(K) \cap \Phi_5 \subset
\partial \Lambda_{3,3}(K) \cap \Phi_6 \equiv
\partial \Lambda_{3,3}(K).
\end{equation}
The first term of this filtration, $\partial \Lambda_{3,3}(K) \cap \Phi_4,$ is
acyclic by Lemma 5. The set $(\partial \Lambda_{3,3}(K) \cap \Phi_5) \sm
(\partial \Lambda_{3,3}(K) \cap \Phi_4)$ is the space of a fiber bundle, whose
base is the space of all generic quadrics in $\CP^2$ (studied in \S~2), and the
fiber over such a quadric $\kappa$ is the open cone
$\stackrel{\circ}{\Lambda}_{3,3}(\kappa).$ The Borel--Moore homology group of
this fiber coincides (up to a shift of dimensions) with the (reduced modulo a
point) homology group of the corresponding link $\partial \Lambda_{3,3}(\kappa)
\sim (S^2)^{*3},$ which is trivial by Lemma 3. Finally, the set $(\partial
\Lambda_{3,3}(K) \cap \Phi_6) \sm (\partial \Lambda_{3,3}(K) \cap \Phi_5)$ is
the space of a fiber bundle, whose base is the space $B(\CP^{2\vee},2)$ of
pairs of complex lines in $\CP^2$, and the fiber over such a pair $\kappa$ is
the open cone $\stackrel{\circ}{\Lambda}_{3,3}(\kappa).$ By the boundary
isomorphism (\ref{bound}), to find the Borel--Moore homology of such a fiber,
we need only to calculate the (usual) homology of the link $\partial
\Lambda_{3,3}(\kappa),$ which already was calculated (see the study of the
column $p=6$) and is trivial.

This proves the triviality of the column $p={\bf 8}$ and terminates the
justification of columns 1 through 10 of Fig.~\ref{sssurf}. The triviality of
the column $p={\bf 11}$ and the degeneration of the spectral sequence at the
term $E^1$ follow from exactly the same reasons as in the proofs of previous
theorems: the former follows from the fact that $\Pi_{3,3} \sm \Sigma_{3,3}$ is
a Stein manifold, and the latter from the fact that it is divisible by $\C^*$
(and from the explicit shape of the term $E^1$).

Proposition 9 and Theorem 4 are completely proved.

\section{Non-degenerate quadratic vector fields in $\C^3$}

{\sc Definition.} A collection of three homogeneous functions
$(v_1,v_2,v_3):\C^3 \to \C^3$ is non-degenerate, if they vanish simultaneously
only at the origin of the issue space.
\medskip

{\sc Theorem 5.} {\em The Poincar\'e polynomial of the rational cohomology
group of the space of non-degenerate triples of homogeneous functions of
degrees 2 in $\C^3$ is equal to $(1+t)(1+t^3)(1+t^5)$.}
\medskip

{\em Proof.} Denote by $\Sigma(2,2,2)$ the discriminant subset in the space
$\C^{18}$ of all triples of quadratic functions in $\C^3$, i.e. the set of all
degenerate triples. The conical resolution of this subset is constructed in
exactly the same way as similar resolutions of spaces $\Sigma_{d,n}$ described
in \S~3; instead of possible singular subsets in $\CP^n$ we use possible sets
of common zeros of a triple of polynomials.
\medskip

{\sc Proposition 10.} {\em There are exactly the following possible sets in
$\CP^2$ of common zeros of triples of homogeneous polynomials of degree 2 in
$\C^3$:

\begin{enumerate}
\item
Any point in $\CP^2$ \quad $\langle 15 \rangle$

\item
Any pair of points \quad $\langle 12 \rangle$

\item
Any generic triple of points (i.e. not lying on the same line) \quad $\langle 9
\rangle$

\item
Any line $\CP^1 \subset \CP^2$ \quad $\langle 9 \rangle$

\item
Any generic quadruple of points \quad $\langle 6 \rangle$

\item
Any line $\CP^1 \subset \CP^2$ plus a point outside it \quad $\langle 6
\rangle$

\item
Any nonsingular quadric in $\CP^2$ \quad $\langle 3 \rangle$

\item
Two lines in $\CP^2$ \quad $\langle 3 \rangle$

\item
Entire $\CP^2$. \quad $\langle 0 \rangle$ \quad $\Box$
\end{enumerate}
}
\medskip

\begin{figure}
\unitlength=1.00mm \thinlines \linethickness{0.4pt}
\begin{picture}(95.00,87.00)
\put(2.00,10.00){\vector(1,0){93.00}} \put(93.00,6.00){\makebox(0,0)[cc]{$p$}}
\put(10.00,4.00){\vector(0,1){83.00}} \put(5.00,85.00){\makebox(0,0)[cc]{$q$}}
\put(18.00,4.00){\line(0,1){77.00}} \put(26.00,81.00){\line(0,-1){77.00}}
\put(34.00,4.00){\line(0,1){77.00}} \put(42.00,81.00){\line(0,-1){77.00}}
\put(50.00,4.00){\line(0,1){77.00}} \put(14.00,73.00){\makebox(0,0)[cc]{${\bf
R}$}} \put(14.00,61.00){\makebox(0,0)[cc]{${\bf R}$}}
\put(14.00,49.00){\makebox(0,0)[cc]{${\bf R}$}}
\put(22.00,49.00){\makebox(0,0)[cc]{${\bf R}$}}
\put(22.00,37.00){\makebox(0,0)[cc]{${\bf R}$}}
\put(22.00,25.00){\makebox(0,0)[cc]{${\bf R}$}}
\put(30.00,13.00){\makebox(0,0)[cc]{${\bf R}$}}
\put(46.00,6.00){\makebox(0,0)[cc]{$5$}}
\put(38.00,6.00){\makebox(0,0)[cc]{$4$}}
\put(30.00,6.00){\makebox(0,0)[cc]{$3$}}
\put(22.00,6.00){\makebox(0,0)[cc]{$2$}}
\put(14.00,6.00){\makebox(0,0)[cc]{$1$}} \put(2.00,16.00){\line(1,0){82.00}}
\put(2.00,22.00){\line(1,0){82.00}} \put(2.00,28.00){\line(1,0){82.00}}
\put(2.00,34.00){\line(1,0){82.00}} \put(2.00,40.00){\line(1,0){82.00}}
\put(2.00,46.00){\line(1,0){82.00}} \put(2.00,52.00){\line(1,0){82.00}}
\put(2.00,58.00){\line(1,0){82.00}} \put(2.00,64.00){\line(1,0){82.00}}
\put(2.00,70.00){\line(1,0){82.00}} \put(2.00,76.00){\line(1,0){82.00}}
\put(5.00,73.00){\makebox(0,0)[cc]{33}} \put(5.00,67.00){\makebox(0,0)[cc]{32}}
\put(5.00,61.00){\makebox(0,0)[cc]{31}} \put(5.00,55.00){\makebox(0,0)[cc]{30}}
\put(5.00,49.00){\makebox(0,0)[cc]{29}} \put(5.00,43.00){\makebox(0,0)[cc]{28}}
\put(5.00,37.00){\makebox(0,0)[cc]{27}} \put(5.00,31.00){\makebox(0,0)[cc]{26}}
\put(5.00,25.00){\makebox(0,0)[cc]{25}} \put(5.00,19.00){\makebox(0,0)[cc]{24}}
\put(5.00,13.00){\makebox(0,0)[cc]{23}} \put(58.00,81.00){\line(0,-1){77.00}}
\put(66.00,4.00){\line(0,1){77.00}} \put(74.00,81.00){\line(0,-1){77.00}}
\put(82.00,4.00){\line(0,1){77.00}} \put(54.00,6.00){\makebox(0,0)[cc]{6}}
\put(62.00,6.00){\makebox(0,0)[cc]{7}} \put(70.00,6.00){\makebox(0,0)[cc]{8}}
\put(78.00,6.00){\makebox(0,0)[cc]{9}}
\end{picture}
\caption{Spectral sequence for non-degenerate quadratic vector fields in
$\C^3$} \label{ssvect}
\end{figure}
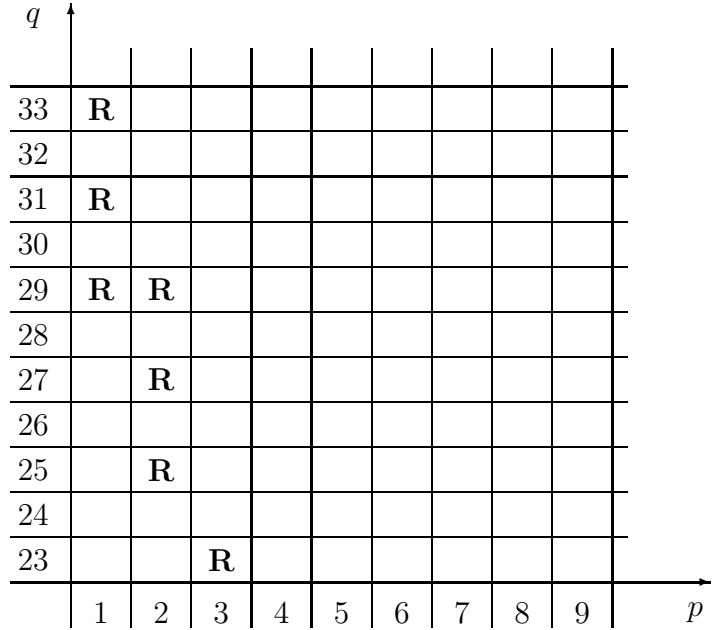

{\sc Proposition 11.} {\em The term $E^1$ of the spectral sequence, generated
by the corresponding filtration and converging to the group $\bar
H_*(\Sigma_(2,2,2) \R),$ is as shown in Fig.~\ref{ssvect}. This spectral
sequence degenerates at this term: $E^\infty \equiv E^1$.}
\medskip

The proof of this proposition repeats essentially some fragments of the proof
of Proposition 7 and will be omitted. The comparison of spectral sequences of
Figs. \ref{ss3} and \ref{ssvect} implies the following version of the
(homological) Smale--Hirsch principle.
\medskip

{\sc Proposition 12.} {\em The gradient embedding $\C^{10} \to \C^{18}$,
mapping any homogeneous polynomial $\C^3 \to \C^1$ of degree 3 to the
collection of its three partial derivatives, induces the isomorphism between
rational cohomology groups of spaces of nonsingular cubical polynomials and
non-degenerate quadratic vector fields in $\C^3$.} \quad $\Box$

\section{Problems and concluding remarks}

1. To explain the 6-dimensional generator of $H^*(N_{4,2})$.
\medskip

2. To find general properties of the above calculations. Of course, these
calculations (and the number of strata to be considered) grow at least
exponentially together with $n$ and $d$. Note however that all the infinite
strata have contributed nothing in our spectral sequences. Is it so also for
greater $n$ and $d$?

Is $H^*(N_{3,n},\R) \simeq H^*(PGL(\CP^n),\R)$ for any $n$?
\medskip

3. To make all the above calculations for the cohomology groups with
coefficients in $\Z$ or at least $\Z_2.$ The same for the real version of the
problem, i.e. for spectral sequences calculating the cohomology of spaces of
nonsingular algebraic hypersurfaces in $\RP^n.$

The involution of complex conjugation acts naturally on our spectral sequences,
thus they are a good tool to check the $M$-property for such spaces or some
strata of the discriminant.

However, the most useful for the original problem of the rigid homotopy
classification should be the similar calculation of $\R$-valued cohomology
groups of spaces of real nonsingular objects. Note that V.~Kharlamov \cite{har}
considered the topology of discriminant sets in a connection with the rigid
classification problem.
\medskip

{\sc Problem}: to express the standard invariants of rigid isotopy
classification of projective curves and surfaces in terms of dual homology
classes of discriminants.
\medskip

4. A collection of points in $\RP^2$ or in $\CP^2$ is called $d$-{\em
sufficient} (respectively, $d$-{\em degenerate}), if any homogeneous polynomial
of degree $d$ with singular points in it vanishes (respectively, has singular
points) at entire curve of smaller degree passing through some of them. The
simplest example: any $[d/2]+1$ (respectively, $d$) points, lying on the same
line, form a $d$-sufficient (respectively, $d$-degenerate) collection.
\medskip

{\sc Problem.} To study the set of such configurations. For which smallest $d$
there are {\it minimal\/} $d$-degenerate configurations such that the
corresponding curve of singularities contains not all their points? For which
$i$ the spaces of all $i$-point $d$-non-degenerate configurations are
disconnected (both in real and complex versions) and how many components do
they have?

To describe the class of algebraic subsets in $\RP^2$ or $\CP^2$ which satisfy
the following $d$-{\em overdeterminacy} property: any homogeneous polynomial of
degree $d$, having singular points at this subset, is identically equal to 0.
\medskip

More generally, let us call a point $y \in \PP^2$ \ a $d$-{\em consequence} of
points $x_1, \ldots, x_l$, if any homogeneous polynomial of degree $d$, having
singularities at these  $l$ points, has also a singularity at the point $y$.
E.g., if $(x_1, \ldots, x_l)$ is a $d$-degenerate (respectively,
$d$-overdetermined) set, then all points of the corresponding curve
(respectively, all points at all) are its $d$-consequences. Are there other
examples?
\medskip

Consider some natural family of algebraic subsets in $\PP^2$, e.g. the set of
all collections, consisting of several curves of fixed degrees and several
points not on these curves. Any point of this family defines a subspace in
$\Pi_d$: the space of polynomials having singularities at all points of the
corresponding algebraic set. Is it possible to describe easily the set of
points of our family, for which the dimension of this subspace ``jumps''?


\begin{thebibliography}{99}

\bibitem{ar70} V.~I.~Arnold, {\it On some topological invariants
of algebraic functions}. Transact. Moscow Math. Society
1970 {\bf 21}, 27--46.


\bibitem{ar89} V.~I.~Arnold, {\it
Spaces of functions with mild singularities.} Funct. Anal. and its Appl.
1989 {\bf 23} (3), 1--10.


\bibitem{aprobl} V.~I.~Arnold, {\em Problems},
Moscow: Phasis, 1999 (in preparation).

\bibitem{borel} A.~Borel. {\it Sur la cohomologie des espaces
fibr\'{e}s principaux et des espaces homog\`{e}nes de groupes de Lie compacts.}
Ann. of Math. 1953 {\bf 57} (2), 115--207.

\bibitem{fuchs} D.~B.~Fuchs. {\it
Cohomology of the braid group
$\mathop{\rm mod} 2$.}
Funct. Anal. and its Appl. 1970 {\bf 4} (2), 62--73.

\bibitem{fuchs2} D.~B.~Fuchs, {\it Maslov--Arnold characteristic
classes}, Doklady (C.R.) of Soviet Ac.Sci., 1968, {\bf 178}, 303--306.

\bibitem{har} V.~M.~Kharlamov, {\it Rigid isotopy classification
of real plane curves of degree 5}, Funct. Anal. and its Appl., 1981, 15:1.

\bibitem{MS} J.~W.~Milnor, J.~D.~Stasheff, {\it Characteristic
Classes,} Princeton, NJ -- Tokyo: Princeton Univ. Press, Univ. of Tokyo Press,
1974, 331 pp. (Ann. of Math. Studies, 76).

\bibitem{viniti} V.A.Vassiliev, {\it Stable cohomology of complements of
discriminants of singularities of smooth functions.
} Current Problems of Math., Newest Results, vol. 33, Itogi Nauki i Tekhniki,
VINITI, Moscow, 1988, 3--29; (Russian, English transl. in  J. Soviet Math. 52:4
(1990), 3217--3230.)

\bibitem{v88} V.~A.~Vassiliev. {\it Lagrange and Legendre
characteristic classes.} 2nd edition. New York a.o.: Gordon \& Breach Science
Publishers, 1993. 265~pp.

\bibitem{v91} V.~A.~Vassiliev, {\it
A geometric realization of the homology of classical Lie  groups,  and
complexes, $S$-dual to the flag manifolds,}
Algebra i Analiz 3:4 (1991), 113--120. (Russian, English transl. in
St.-Petersburg Math. J., 3:4 (1992), 809--815.)

\bibitem{book} V.~A.~Vassiliev, {\it Complements of discriminants of smooth
maps: topology and applications. Revised ed.}, Providence RI: AMS,  1994.
(Transl. of Math. Monographs, 98)

\bibitem{congress} V.~A.~Vassiliev, {\it Topology of discriminants
and their complements,} Proceedings of the Intern. Congress of Math. (Z\"urich,
1994); Birkh\"auser Verlag, Basel; 1995, 209--226.

\bibitem{phasis} V.~A.~Vassiliev, {\em Topology of complements of
discriminants,} Moscow: Phasis, 1997 (in Russian).


\end{thebibliography}
\end{document}